\newcommand{\pbaddress}{biran@math.tau.ac.il}

\documentclass[12pt]{amsart}

\usepackage{epsfig}
\usepackage{amscd}
\usepackage[mathscr]{eucal}
\usepackage{amssymb}
\usepackage{amsxtra}
\usepackage{amsmath}
\usepackage[all]{xy}

\theoremstyle{plain}
\newtheorem{mainthm}{Theorem}
\newtheorem{thm}{Theorem}[subsection]

\newtheorem{lem}[thm]{Lemma}
\newtheorem{prop}[thm]{Proposition}

\theoremstyle{definition}

\theoremstyle{remark}
\newtheorem{rem}[thm]{Remark}
\newtheorem*{remnonum}{Remark}

\newtheorem*{remsnonum}{Remarks}

\newtheorem*{exnonum}{Example}
\newtheorem*{exsnonum}{Examples}

\newtheorem*{qssnonum}{Questions}

\newcommand{\cntrs}{\setcounter{thm}{0}
  \renewcommand{\thethm}{\thesection.\Alph{thm}}}
\newcommand{\cntrsb}{\setcounter{thm}{0}
  \renewcommand{\thethm}{\thesubsection.\Alph{thm}}}

\newcommand{\Id}{{{\mathchoice {\rm 1\mskip-4mu l} {\rm 1\mskip-4mu l}
      {\rm 1\mskip-4.5mu l} {\rm 1\mskip-5mu l}}}}

\begin{document}

\title{Lagrangian non-Intersections} 

\date{\today}

\thanks{This research was supported by the ISRAEL SCIENCE FOUNDATION
  (grant No.  205/02 *)}

\author{Paul Biran} \address{Paul Biran, School of Mathematical
  Sciences, Tel-Aviv University, Ramat-Aviv, Tel-Aviv 69978, Israel}
\email{\pbaddress}


%
%
% ...

\maketitle

%
%

\section{Introduction and main results} \cntrs
\label{S:Intro}
The present paper is devoted to Lagrangian submanifolds of symplectic
manifolds and their intersection patterns.

One of the cornerstones of symplectic topology is the rigidity of
symplectic structures reflected in the behavior of Lagrangian
submanifolds and their intersections. The first non-trivial
restrictions on Lagrangian embeddings go back to Gromov's work
of~1985~\cite{Gr:phol}.  One of the many results of that paper is the
following fundamental theorem: {\sl In $\mathbb{R}^{2n}$ there are no
  closed exact Lagrangian submanifolds, in particular there are no
  closed Lagrangian submanifolds $L \subset \mathbb{R}^{2n}$ with
  $H^1(L;\mathbb{R})=0$.}

Since then, Gromov's techniques have been extended and, in combination
with other methods, have led to more restrictions on the topology of
Lagrangian submanifolds, most of the results being for Lagrangians
submanifolds of $\mathbb{R}^{2n}$ and of cotangent bundles (see
e.g.~\cite{ALP,Bu:TSn,El:symp90,La-Si:Lag,Moh:lk,Moh:Kbottle,Oh:spectral,
  Oh:relative,Po:Maslov,Se:TSn,Vi:Maslov,Vi:functors1} for a partial
list of older and newer results).

Only relatively recently first results on the topology of Lagrangians
in closed manifolds have been obtained by Seidel~\cite{Se:graded} and
later on by Biran and Cieliebak~\cite{Bi-Ci:closed}. Note that when
studying Lagrangians in an arbitrary manifold one encounters all
Lagrangian submanifolds of $\mathbb{R}^{2n}$. This is due to the fact
that every symplectic manifold $M^{2n}$ is locally modeled on
$\mathbb{R}^{2n}$, hence every Lagrangian submanifold $L \subset
\mathbb{R}^{2n}$ can also be Lagrangianly embedded into $M^{2n}$.
Thus, Lagrangian embeddings into $\mathbb{R}^{2n}$ should, in a sense,
be regarded as the local case. However, our understanding of
Lagrangian submanifolds of $\mathbb{R}^{2n}$ is still quite limited,
in particular also that of ``local'' Lagrangian submanifolds in any
symplectic manifold $M^{2n}$. (Thus, in this case, ``local'' turns out
to be difficult.)

In this paper we concentrate on ``global'' Lagrangian submanifolds.
One way to ``mod out'' the local Lagrangians is to assume for example
that the first homology of the Lagrangians is either zero or torsion.
In view of the preceding theorem of Gromov such Lagrangian
submanifolds must be global in the sense that they cannot lie entirely
in a Darboux chart.

One of the phenomena arising from our results below is that in certain
symplectic manifolds the topology of Lagrangians with small first
homology is extremely restricted. It turns out that in some cases
(e.g. $M={\mathbb{C}}P^n$) certain assumptions on $H_1(L;\mathbb{Z})$
of a Lagrangian $L$ completely determine the entire homology of $L$.
This phenomenon is illustrated in Theorems~\ref{T:cpn}-~\ref{T:cpncpn}
of Section~\ref{Sb:homological} below. Note that very recently
examples of this phenomenon have been discovered also for cotangent
bundles of spheres by Buhovski~\cite{Bu:TSn} and by
Seidel~\cite{Se:TSn}.

The second phenomenon presented in this paper belongs to the framework
of Lagrangian intersections. Our results show that certain symplectic
manifolds $M$ contain some kind of ``Lagrangian core'' $\Lambda
\subset M$ which dominates intersections, in the sense that many
global Lagrangians $L \subset M$ must intersect $\Lambda$, the
intersection points being irremovable by symplectic diffeomorphisms.
Results in this direction are presented in Section~\ref{Sb:intersect}.

\subsection{Homological uniqueness of Lagrangian submanifolds} \cntrsb 
\label{Sb:homological}
Here and in what follows all Lagrangian submanifolds are assumed to be
compact and without boundary, unless otherwise explicitly stated.

\subsubsection*{Lagrangian submanifolds of ${\mathbb{C}}P^n$}
Let ${\mathbb{C}}P^n$ be the complex projective space, endowed with
its standard K\"{a}hler structure. It is well known that
${\mathbb{C}}P^n$ has no Lagrangian submanifolds $L$ with
$H_1(L;\mathbb{Z})=0$ (see~\cite{Se:graded}, see
also~\cite{Bi-Ci:closed}). Note however that there do exist
Lagrangians $L \subset {\mathbb{C}}P^n$ with $H_1(L; \mathbb{Z})$
torsion. For example, the real projective space
$$\mathbb{R}P^n \approx \big\{[z_0:\ldots:z_n] \in {\mathbb{C}}P^n
\mid z_i \in \mathbb{R} \; \forall i\, \big\} \subset
{\mathbb{C}}P^n$$
is such a Lagrangian submanifold (for $n\geq 2$). In
fact Seidel proved in~\cite{Se:graded} that {\sl every Lagrangian
  submanifold $L \subset {\mathbb{C}}P^n$ with
  $H^1(L;\mathbb{Z}_{2n+2})=\mathbb{Z}_2$ must satisfy
  $H^*(L;\mathbb{Z}_2) \cong H^*(\mathbb{R}P^n;\mathbb{Z}_2)$ as
  graded vector spaces}. Below we shall prove a stronger statement
which gives information also on the cohomology {\em ring} of $L$.
Henceforth we say that an abelian group $H$ is $q$-torsion if for
every $\alpha \in H$ we have $q\alpha=0$. (This, by our conventions,
{\em includes} the case when $H$ is the trivial group.) Our first
result is:
\begin{mainthm} \label{T:cpn}
   Let $L \subset {\mathbb{C}}P^n$ be a Lagrangian submanifold with
   $H_1(L;\mathbb{Z})$ $2$-torsion. Then:
   \begin{enumerate}
     \item There exists an isomorphism of graded vector spaces
      $H^*(L;\mathbb{Z}_2) \cong H^*(\mathbb{R}P^n; \mathbb{Z}_2)$.
      Moreover, if $a \in H^2({\mathbb{C}}P^n; \mathbb{Z}_2)$ is the
      generator then $a|_{_L} \in H^2(L;\mathbb{Z}_2)$ generates the
      subalgebra $H^{\textnormal{even}}(L;\mathbb{Z}_2)$, and
      $H^{\textnormal{odd}}(L; \mathbb{Z}_2) = H^1(L; \mathbb{Z}_2)
      \cup H^{\textnormal{even}}(L; \mathbb{Z}_2)$.
     \item When $n$ is even, the isomorphism $H^*(L;
      \mathbb{Z}_2)\cong H^*(\mathbb{R}P^n; \mathbb{Z}_2)$ is in fact
      an isomorphism of graded algebras.
   \end{enumerate}
\end{mainthm}
Other than $\mathbb{R}P^n \subset {\mathbb{C}}P^{n+1}$ we are not
aware of any example of Lagrangians $L$ with $2$-torsion
$H_1(L;\mathbb{Z})$. However, Chiang~\cite{Chiang:RP3} constructed an
example of a Lagrangian in $L \subset {\mathbb{C}}P^3$ with
$H_i(L;\mathbb{Z}_2)=\mathbb{Z}_2$ for every $i$. This Lagrangian is a
quotient of $\mathbb{R}P^3$ by the dihedral group $D_3$. It satisfies
$H_1(L;\mathbb{Z})=\mathbb{Z}_4$.  See~\cite{Chiang:RP3} for the
details as well as other interesting examples of Lagrangian
submanifolds of complex projective spaces of various dimensions.

\subsubsection*{Lagrangian submanifolds of ${\mathbb{C}}P^n\times X$}
It has been proved in~\cite{Bi-Ci:closed} that if $X$ is a closed
symplectic manifold with $\pi_2(X)=0$ then for $n\geq
\dim_{\mathbb{C}}X$, ${\mathbb{C}}P^n\times X$ has no simply connected
Lagrangian submanifolds. On the other hand, if $n<\dim_{\mathbb{C}}X$,
then ${\mathbb{C}}P^n\times X$ may have such Lagrangians. Indeed, for
any symplectic manifold $X$ of $\dim_{\mathbb{C}}X=n+1$,
${\mathbb{C}}P^n \times X$ has a Lagrangian sphere, after a possible
rescaling of the symplectic form on the $X$ factor
(see~\cite{Bi-Ci:closed} for details). The following theorem shows
that homologically this is the only example.
\begin{mainthm} \label{T:cpnX}
   Let $X$ be a symplectic manifold of $\dim_{\mathbb{C}}X=n+1$ with
   $\pi_2(X)=0$ and which is either closed or has a symplectically
   convex end. If $L^{2n+1}\subset {\mathbb{C}}P^n\times X$, $n\geq
   1$, is a simply connected Lagrangian submanifold then
   $H^*(L^{2n+1}; \mathbb{Z}_2) \cong H^*(S^{2n+1}; \mathbb{Z}_2)$.
\end{mainthm}

In contrast to the result of~\cite{Bi-Ci:closed} mentioned above, if
one drops the condition ``$\pi_2(X)=0$'', then ${\mathbb{C}}P^n\times
X$ may have simply connected Lagrangian submanifolds even if $n\geq
\dim_{\mathbb{C}} X$. For example take $X$ to be ${\mathbb{C}}P^n$ and
consider ${\mathbb{C}}P^n \times {\mathbb{C}}P^n$ endowed with the
equally weighted split standard symplectic structure. Then
$L={\mathbb{C}}P^n$ embeds Lagrangianly as the ``anti-diagonal'' in
${\mathbb{C}}P^n\times {\mathbb{C}}P^n$, namely
$${\mathbb{C}}P^n \ni [z_0: \ldots: z_n] \longmapsto \big([z_0:
\ldots: z_n], [\overline{z}_0: \ldots: \overline{z}_n]\big) \in
{\mathbb{C}}P^n\times {\mathbb{C}}P^n.$$
The following theorem shows
that, again, homologically this is the only example.
\begin{mainthm} \label{T:cpncpn}
   Let $L \subset {\mathbb{C}}P^n \times {\mathbb{C}}P^n$ be a
   Lagrangian submanifold with $H_1(L; \mathbb{Z})=0$. Then
   $H^*(L;\mathbb{Z}_2) \cong H^*({\mathbb{C}}P^n; \mathbb{Z}_2)$ as
   graded algebras. Moreover, if $a\in H^2({\mathbb{C}}P^n\times
   {\mathbb{C}}P^n; \mathbb{Z}_2)$ is the generator of $H^2$ of either
   factor of ${\mathbb{C}}P^n\times {\mathbb{C}}P^n$ then $a|_{_L} \in
   H^2(L;\mathbb{Z}_2)$ generates $H^*(L; \mathbb{Z}_2)$ as an
   algebra.
\end{mainthm}

Again, besides $L={\mathbb{C}}P^n$ we are not aware of any other
examples of Lagrangian submanifolds in ${\mathbb{C}}P^n \times
{\mathbb{C}}P^n$ with $H_1(L;\mathbb{Z})=0$.

\subsection{Lagrangian spheres} \cntrsb \label{Sb:spheres}
Here we present new restrictions on Lagrangian embeddings of spheres.
Let us mention that Lagrangian spheres appear in various mathematical
contexts other than symplectic geometry (e.g.  singularity
theory~\cite{Bi:symp-obst}) and thus deserve special attention beyond
the scope of symplectic geometry.

In view of the fact that manifolds of the type ${\mathbb{C}}P^n\times
X$ with $\dim_{\mathbb{C}}X=n+1$ have Lagrangian spheres it makes
sense to ask what happens for $X$ of other dimensions. The next
theorem gives a partial answer to this question.
\begin{mainthm} \label{T:cpnX-sphere}
   Let $X$ be a symplectic manifold with $\pi_2(X)=0$ which is either
   closed or has a symplectically convex end. If
   ${\mathbb{C}}P^n\times X$ has a Lagrangian sphere (where $n,\dim
   X>0$) then $\dim_{\mathbb{C}}X \equiv n+1(\bmod{2n+2})$.
\end{mainthm}

Recall that a symplectic manifold $(M, \omega)$ is called {\em
  spherically monotone} if $[\omega]|_{\pi_2(M)} \neq 0$ and there
exists $\lambda > 0$ such that $[\omega]=\lambda c_1^M$ on $\pi_2(M)$.
Here, and in what follows, $c_1^M$ stands for the first Chern class of
the tangent bundle of $M$, viewed (in a canonical way) as a complex
vector bundle. We denote by $N_M \in \mathbb{Z}_+$ the {\em minimal
  Chern number}, namely $N_M=\min \{ c_1^M(A) \mid A\in \pi_2(M),
c_1^M(A)>0 \}$.

\begin{mainthm} \label{T:M-cover}
   Let $X$ be a symplectic manifold that can be covered by a
   symplectic manifold which is symplectomorphic to a domain in
   $\mathbb{C}^m$.  Let $M$ be a spherically monotone closed
   symplectic manifold.  Assume $\dim X, \dim M >0$. If $M \times X$
   has a Lagrangian sphere then $2N_M \mid \dim_{\mathbb{C}}M +
   \dim_{\mathbb{C}}X + 1$.
\end{mainthm}
Examples of manifolds $X$ satisfying the conditions of
Theorem~\ref{T:M-cover} include symplectic tori, and ball quotients,
both endowed with K\"{a}hler symplectic structures.
Theorem~\ref{T:M-cover} is in fact a special case of the more general
Theorem~\ref{T:M-cover-2} which will be proved in
Section~\ref{S:Generalizations}.

\begin{mainthm} \label{T:cpnM-sphere}
   Let $(M, \omega)$ be a spherically monotone closed symplectic
   manifold with $[\omega] = c_1^M$ on $\pi_2(M)$ and denote
   $m=\dim_{\mathbb{C}} M$. Let ${\mathbb{C}}P^n\times M$ be endowed
   with the symplectic form $\Omega=(n+1)\sigma \oplus \omega$, where
   $\sigma$ is the standard symplectic K\"{a}hler form of
   ${\mathbb{C}}P^n$ normalized so that the area of a projective line
   is $1$. Suppose that $({\mathbb{C}}P^n\times M, \Omega)$ has a
   Lagrangian sphere, where $n+m\geq 3$. Then $2\gcd{(n+1,N_M)}\mid
   n+m+1$.
\end{mainthm}
\begin{exnonum}
   Theorem~\ref{T:cpnM-sphere} implies that if $({\mathbb{C}}P^n
   \times {\mathbb{C}}P^m, (n+1)\sigma \oplus (m+1)\sigma)$ has a
   Lagrangian sphere for $n+m\geq 3$ then $\gcd{(n+1,m+1)}=1$ and
   $n+m$ is odd.
\end{exnonum}

\begin{remnonum}
   Embarrassingly, the only example known to us of Lagrangian spheres
   in manifolds of the types appearing in
   Theorems~\ref{T:cpnX-sphere}-\ref{T:cpnM-sphere}, are all in
   manifolds of the type ${\mathbb{C}}P^n \times X$ where
   $\dim_{\mathbb{C}} X=n+1$. It would be interesting to figure out
   for example whether or not ${\mathbb{C}}P^n \times X$ admits a
   Lagrangian sphere (or even a homology sphere) when
   $\dim_{\mathbb{C}}X=(2k+1)(n+1)$, $k \geq 1$.
\end{remnonum}

\subsection{Lagrangian intersections} \cntrsb \label{Sb:intersect}

Here we present new results on Lagrangian intersections.  The pattern
that stands out in all the examples below is the existence of a
``core'' $\Lambda$ consisting of a finite union of (possibly) immersed
Lagrangian spheres with the property that every Lagrangian submanifold
with prescribed topological properties (e.g. simply connected) must
intersect $\Lambda$.

\subsubsection*{Intersections in the quadric} 
Let $Q^n \subset {\mathbb{C}}P^{n+1}$ be the complex $n$-dimensional
smooth quadric defined by the equation:
$$Q^n = \left\{ [z_0: \ldots : z_{n+1}] \in {\mathbb{C}}P^{n+1} \;
   \biggm| z_0^2 + \cdots + z_{n+1}^2 =0 \right\},$$
endowed with the
symplectic structure induced from ${\mathbb{C}}P^{n+1}$. Let
$\Lambda_Q \subset Q^n$ be the corresponding ``real'' quadric, namely:
$$\Lambda_Q = \left\{[z_0: \ldots : z_{n+1}] \in Q^n \biggm| z_0 \in
   \mathbb{R}, z_1, \ldots, z_{n+1} \in i\mathbb{R} \right\}.$$
It is
not hard to see that $\Lambda_Q \subset Q^n$ is in fact a Lagrangian
sphere.

Let $A_{Q}^* = \oplus_{i=0}^n A_Q^i$ be the following graded vector
space over $\mathbb{Z}_2$:
\begin{equation*}
   \begin{cases}
      A_Q^0=\mathbb{Z}_2, \; A_Q^1=\mathbb{Z}_2, \\
      A_Q^i=0 \quad \textnormal{for every } 1<i<n-1, \\
      A_Q^{n-1}=\mathbb{Z}_2, \; A_Q^n=\mathbb{Z}_2.
   \end{cases}
\end{equation*}

\begin{mainthm} \label{T:AQ}
   Let $L\subset Q^n$, $n\geq 3$, be a Lagrangian submanifold such
   that $H_1(L;\mathbb{Z})$ is a $2$-torsion group.  If
   $H^*(L;\mathbb{Z}_2)$ is not isomorphic to $A_Q^*$ then $L \cap
   \Lambda_Q \neq \emptyset$. In particular for every Lagrangian
   submanifold $L \subset Q^n$ with $H_1(L;\mathbb{Z})=0$, we have $L
   \cap \Lambda_Q \neq \emptyset$.
\end{mainthm}

\begin{remsnonum}
   \begin{enumerate}
     \item Note that the intersection between $L$ and $\Lambda_Q$ is
      in general not due to topological reasons but rather to
      symplectic ones. For example when $n$ is odd, every Lagrangian
      sphere (in particular $\Lambda_Q$) can be displaced from itself
      via an arbitrary small (non symplectic) diffeomorphism.
     \item The quadric $Q^n$ has many different Lagrangian
      submanifolds with $H_1(L; \mathbb{Z})$ either zero or
      $2$-torsion.  First of all it has simply connected Lagrangians
      (e.g.  $\Lambda_Q$ itself).  Next, for every $0\leq r \leq n$
      consider
      $$L_r = \left\{ [z_0: \ldots: z_{n+1}] \in Q^n \; \big|\; z_0,
         \ldots, z_r \in \mathbb{R}, \quad z_{r+1}, \ldots, z_{n+1}
         \in i \mathbb{R} \right\}.$$
      It is not hard to see that all
      the $L_r$'s are Lagrangian submanifolds of $Q^n$ and that $L_r$
      is diffeomorphic to $(S^r \times S^{n-r})/ \mathbb{Z}_2$, where
      $\mathbb{Z}_2$ acts on both factors by the antipode map. It
      easily follows that $L_0$ and $L_n$ are spheres, and a simple
      computation shows that for every $1< r < n-1$, we have
      $H_1(L_r;\mathbb{Z}) = \mathbb{Z}_2$ but $H^2(L_r; \mathbb{Z}_2)
      \neq 0$. Therefore when $n > 3$, $H^*(L_r; \mathbb{Z}_2)$ is not
      isomorphic to $A_Q^*$. It follows from Theorem~\ref{T:AQ} that
      for $n > 3$ any Lagrangian submanifold $L \subset Q^n$
      diffeomorphic to one of the $L_r$'s must intersect $\Lambda_Q$.
      In particular the intersection $L_r \cap \Lambda_Q \neq
      \emptyset$ cannot be removed via a symplectic isotopy.
     \item We do not know of any examples of Lagrangians $L$ with
      $2$-torsion $H_1(L;\mathbb{Z})$ lying in the complement of
      $\Lambda_Q$.
     \item The statement of Theorem~\ref{T:AQ} remains true if instead
      of assuming $H_1(L;\mathbb{Z})$ is $2$-torsion one assumes that
      $L$ is monotone with minimal Maslov number $N_L=n$. (See the
      second remark at the end of the proof of the Theorem in
      Section~\ref{S:Proofs-main}.) It would be interesting to figure
      out if at least there are any monotone Lagrangians $L \subset
      Q^n$ (not necessarily with $2$-torsion $H_1$) with $N_L=n$ that
      lie in complement of $\Lambda_Q$.
   \end{enumerate}
\end{remsnonum}

\begin{qssnonum}
   \begin{enumerate}
     \item Theorem~\ref{T:AQ} implies that (for $n \geq 3$) every
      Lagrangian sphere $L \subset Q^n$, must intersect $\Lambda_Q$.
      Is it true that every two Lagrangian spheres $L_1, L_2 \subset
      Q^n$ must intersect each other? Note that when $n=$ even this
      easily follows by computing the intersection number $[L_1] \cdot
      [L_2]$, but for $n=$ odd this does not seem to follow from
      purely topological reasons. An affirmative answer would have the
      following consequences in algebraic geometry: the complex
      quadric (of dimension $\geq 3$) cannot be degenerated to a
      variety having two or more isolated singularities.
      See~\cite{Bi:ICM2002, Bi:ECM2004, Bi:symp-obst} for more details
      and precise statements.
     \item Note that the cohomology of the Lagrangian $L_1 \subset
      Q^n$ (taken from (2) above with $r=1$) is precisely $A_Q^*$. Can
      $L_1$ be Hamiltonianly isotoped to lie in the complement of
      $\Lambda_Q$ ?  It is known that $L_1$ can also be embedded as a
      monotone Lagrangian in $\mathbb{C}^n$ with minimal Maslov number
      $=n$ (see~\cite{Po:Maslov}). Can $L_1$ be embedded as a {\em
        monotone} Lagrangian in $Q^n \setminus \Lambda_Q$ ?
     \item Theorem~\ref{T:AQ} will not be proved by showing that the
      Floer homology $HF(L,\Lambda_Q)$ is not zero. It would we be
      interesting to figure out whether this is indeed so.
   \end{enumerate}
\end{qssnonum}

\subsubsection*{Intersections in hypersurfaces of ${\mathbb{C}}P^{n+1}$}
Theorem~\ref{T:AQ} can be generalized as follows. Let $\Sigma_d^{n}
\subset {\mathbb{C}}P^{n+1}$ be a smooth complex hypersurface of
degree $d$, viewed as a symplectic manifold endowed with the
symplectic structure induced from ${\mathbb{C}}P^{n+1}$.

\begin{mainthm} \label{T:hypersurface}
   For every $d>2$ there exist $d^{n+1}$ (possibly) immersed
   Lagrangian spheres $S_1, \ldots, S_{d^{n+1}} \subset \Sigma_d^{n}$
   such that their union $\Lambda_d=S_1 \cup \ldots \cup S_{d^{n+1}}$
   has the following properties:
   \begin{enumerate}
     \item When $d \leq \frac{n+1}{2}$ or $d \geq \frac{3}{2}(n+1)$,
      $n\geq 3$, every Lagrangian submanifold $L\subset \Sigma_d^n$
      with $H_1(L;\mathbb{Z})=0$ must satisfy $L \cap \Lambda_d \neq
      \emptyset$.
     \item When $d \leq \frac{n+1}{2}$ and $n=$ even, for every
      Lagrangian submanifold $L \subset \Sigma_d^n$ with
      $H_1(L;\mathbb{Z})$ $2$-torsion either $L \cap \Lambda_d \neq
      \emptyset$ or $L$ has the following properties:
      \begin{enumerate}
        \item $H^d(L;\mathbb{Z}_2)= \ldots =
         H^{n-d}(L;\mathbb{Z}_2)=0$.
        \item $\beta_i(L) = \beta_{d-1-i}(L) = \beta_{i+1+n-d}(L) =
         \beta_{n-i}(L)$ for every $0\leq i \leq d-1$, where
         $\beta_j(L)$ stands for the $j$'th $\mathbb{Z}_2$-Betti
         number of $L$.
      \end{enumerate}
      
     \item When $d \leq \frac{n+1}{2}$ and $d=$ odd, for every
      Lagrangian submanifold $L \subset \Sigma_d^n$ with
      $H_1(L;\mathbb{Z})$ $2$-torsion either $L \cap \Lambda_d \neq
      \emptyset$ or $L$ has the following properties:
      \begin{enumerate}
        \item $H^k(L;\mathbb{Z}_2) \xrightarrow{\cup w}
         H^{k+2}(L;\mathbb{Z}_2)$ is an isomorphism for every $d\leq k
         \leq n-d-2$.
        \item $H^{d-1}(L;\mathbb{Z}_2) \xrightarrow{\cup w}
         H^{d+1}(L;\mathbb{Z}_2)$ is surjective.
         
         Here, $w \in H^2(L;\mathbb{Z}_2)$ is the restriction of the
         generator $a \in H^2({\mathbb{C}}P^{n+1}; \mathbb{Z}_2)$ to
         $L \subset \Sigma_d^n \subset {\mathbb{C}}P^{n+1}$.
      \end{enumerate}
     \item When $2 < d \leq n+1$, $n \geq 3$ and $2(n+2-d) \nmid n+1$,
      every Lagrangian sphere $L \subset \Sigma_d^n$ must satisfy $L
      \cap \Lambda_d \neq \emptyset$.
     \item Let $d,t \geq 2$, $n\geq 3$ be such that $d \geq
      \frac{t(n-1)}{2} + n+2$. Then every Lagrangian submanifold with
      $H_1(L;\mathbb{Z})$ $t$-torsion must satisfy $L \cap \Lambda_d
      \neq \emptyset$.
   \end{enumerate}

\end{mainthm}
We have not been able to explicitly compute the Lagrangian spheres in
$\Lambda_d$, however in Section~\ref{sbsb:deg-d} below we show that
$S_2, \ldots, S_{d^{n+1}}$ are all obtained from $S_1$ by applying
suitable automorphisms of $\Sigma_d^n$.

Examples of Lagrangian submanifolds $L \subset \Sigma_d^n$ that
satisfy the conditions of the Theorem~\ref{T:hypersurface} come from
Picard-Lefschetz theory. Indeed $\Sigma_d^n$, $d\geq 2$, can be
included as a fibre in a degeneration with isolated singularities,
hence by the Lagrangian vanishing cycle construction it must contain
Lagrangian spheres
(see~\cite{Ar:monodromy,Do:polynom,Kh-Se:quivers,Se:vcycles-mut,Se:thesis},
see also~\cite{Bi:symp-obst}). Here are more explicit examples of
Lagrangians in $\Sigma_d^n$: write $\Sigma_d^n$ as $\left\{ z_0^d +
   \cdots + z_{n+1}^d = 0 \right\} \subset {\mathbb{C}}P^{n+1}$. Let
$\tau \in \mathbb{C}$ be a root of $-1$ of order $d$.  Then it is easy
to see that when $d=$ even
$$\left\{ [z_0: \ldots: z_{n+1}] \in \Sigma_d^n \biggm| z_0 \in \tau
   \mathbb{R}, z_j \in \mathbb{R} \quad \textnormal{for every } 1\leq
   j \leq n+1 \right\}$$
is a Lagrangian sphere. When $d=$ odd we also
have Lagrangians homeomorphic to $\mathbb{R}P^n$, namely $\Sigma_d^n
\cap \mathbb{R}P^{n+1}$. (See Appendix~A of~\cite{Laz:discs} for an
explicit homeomorphism.)

\begin{qssnonum}
   \begin{enumerate}
     \item Theorem~\ref{T:hypersurface} implies in particular that
      (for some values of $d$) it is impossible to embed a Lagrangian
      sphere which is disjoint from the spheres $S_1, \ldots,
      S_{d^{n+1}}$. In view of this one is led to speculate that the
      maximal number of mutually disjoint Lagrangian spheres in
      $\Sigma_d^n$ is finite, or even that this number is not bigger
      than $d^{n+1}$ (c.f. the first question after
      Theorem~\ref{T:AQ}.)
     \item Under the conditions of statement~(2) of
      Theorem~\ref{T:hypersurface}, Lagrangians $L \subset \Sigma_d^n
      \setminus \Lambda_d$ with $2$-torsion $H_1(L;\mathbb{Z})$ must
      have the same $\mathbb{Z}_2$-coefficients cohomology as a
      manifold of the type $L_0 \times S^{n+1-d}$, where $L_0$ is a
      $(d-1)$-dimensional manifold. Can such manifolds be Lagrangianly
      embedded in $\Sigma_d^n$ ? In $\Sigma_d^n \setminus \Lambda_d$ ?
     \item Consider the manifold $L=(S^1 \times
      S^{n+1-d})/\mathbb{Z}_2 \times S^{d-2}$, where $\mathbb{Z}_2$
      acts on both factors of $S^1 \times S^{n+1-d}$ by the antipode
      map. Note that $L$ also satisfies the cohomological restrictions
      predicted by statement~(2) of Theorem~\ref{T:hypersurface}. It
      is known that $L$ admits a monotone Lagrangian embedding into
      $\mathbb{C}^n$ with minimal Maslov number $=n+2-d$
      (see~\cite{Po:Maslov}). Does $L$ admit a {\em monotone}
      Lagrangian embedding into $\Sigma_d^n$ ? Into $\Sigma_d^n
      \setminus \Lambda_d$~?
   \end{enumerate}
\end{qssnonum}

\subsubsection*{Intersections in a hypersurface of 
  ${\mathbb{C}}P^n\times {\mathbb{C}}P^n$}

Consider ${\mathbb{C}}P^n\times {\mathbb{C}}P^n$ endowed with the
standard symplectic form $\sigma_{\textnormal{std}} \oplus
\sigma_{\textnormal{std}}$ and let $\Sigma^{2n-1} \subset
{\mathbb{C}}P^n \times {\mathbb{C}}P^n$ be the complex hypersurface
defined by the equation
$$\Sigma^{2n-1} = \left\{ \sum_{j=0}^{n-1} z_j w_j = z_n w_n
\right\},$$
where $[z_0: \ldots :z_n]$, $[w_0: \ldots :w_n]$ are
homogeneous coordinates on ${\mathbb{C}}P^n \times {\mathbb{C}}P^n$.
We endow $\Sigma^{2n-1}$ with the symplectic structure induced from
${\mathbb{C}}P^n\times {\mathbb{C}}P^n$. Put
$$\Lambda_{\Sigma} = \Big\{ ([z_0: \ldots :z_n], [w_0: \ldots :w_n])
\in \Sigma^{2n-1} \bigm| w_j \overline{z}_n = \overline{z}_j w_n \;
\forall \; 0\leq j \leq n-1 \Big\}.$$
A simple computation shows that
$\Lambda_{\Sigma} \subset \Sigma^{2n-1}$ is a Lagrangian sphere.

\begin{mainthm} \label{T:sigmacpncpn}
   Let $L^{2n-1} \subset \Sigma^{2n-1}$, $n\geq 2$, be a Lagrangian
   submanifold with $H_1(L; \mathbb{Z})=0$. Then either
   $H^*(L;\mathbb{Z}_2) \cong H^*(S^{2n-1}; \mathbb{Z}_2)$ or $L\cap
   \Lambda_{\Sigma} \neq \emptyset$.
\end{mainthm}

\subsection{Discussion} \label{S:discussion}

The phenomenon arising in the results of Section~\ref{Sb:intersect}
above is the existence of certain ``Lagrangian subsets'' that dominate
intersections in the sense that {\em all} Lagrangian submanifolds with
specified topology must intersect them.  A similar phenomenon is
already known for cotangent bundles, where a result due to
Gromov~\cite{Gr:phol} implies that every Lagrangian submanifold of a
cotangent bundle $L \subset T^*(X)$ with $H^1(L;\mathbb{R})=0$ must
intersect the zero section. Our results show that such a phenomenon
holds also in {\em closed} manifolds.  It is interesting to note that,
similarly to our case, also for cotangent bundles it is currently
unknown whether or not the Floer homology $HF(L,O_X)$ of $L$ and the
zero section $O_X$ is non-trivial.  (Compare Question~(2) after
Theorem~\ref{T:AQ}.)

Finally, let us remark that the above intersection phenomena, in
general, seize to hold in the $C^{\infty}$-category. Indeed in each of
the Theorems~G-I one can usually remove all the intersection points by
a smooth diffeomorphism.

\subsection{Methods and ideas}
The methods and tools used in this paper consist of two main
ingredients. The first one is Floer theory for Lagrangian
submanifolds. In particular we use the extension of Floer homology to
monotone Lagrangian submanifolds due to
Oh~\cite{Oh:HF1,Oh:spectral,Oh:relative} which gives rise to an an
algebraic approach for computing Floer homology in terms of a spectral
sequence.

The second ingredient is a geometric technique, developed by the
author in~\cite{Bi:barriers} and in this paper, by which it is
sometimes possible to compute Floer homology in a geometric way. Our
techniques enable to perform certain transformations to a given
Lagrangian submanifold resulting in a new Lagrangian that can be
Hamiltonianly displaced. In particular, we obtain vanishing of Floer
homology.  This vanishing combined with the algebraic computations
mentioned above is the key point behind most of our results.

Let us now describe in some more detail the main ideas of the paper.
Let $\Sigma^{2n}$ be a closed symplectic manifold. We shall
concentrate on the case when $\Sigma$ can be symplectically embedded
as a hyperplane section of a higher dimensional symplectic manifold,
say $M^{2n+2}$. The idea is that in view of the decomposition
technique developed by the author in~\cite{Bi:barriers} the symplectic
properties of $M \setminus \Sigma$ can be used to study the symplectic
topology of $\Sigma$ itself. Consider a tubular neighbourhood $U$ of
$\Sigma$ in $M$. Its boundary $\partial U$ is a circle bundle $\pi:
\partial U \to \Sigma$ over $\Sigma$. Now let $L \subset \Sigma$ be a
Lagrangian submanifold.  As we shall see in
Section~\ref{Sb:Lag-circle} below, if we chose $U$ carefully the
restriction of this circle bundle to $L$, is a Lagrangian submanifold
$\Gamma_L = \pi^{-1}(L)$ lying in $M \setminus \Sigma$. Note that $M
\setminus \Sigma$ is a Stein manifold.

The symplectic theory of Stein manifolds now comes into play. Recall
that Stein manifolds are divided into subcritcal and critical ones
(see Section~\ref{S:Stein} for the precise definitions). Subcritical
Stein manifolds have the feature that (after a suitable completion)
all compact subsets can be Hamiltonianly displaced. In particular the
Floer homology of Lagrangian submanifolds must vanish.

Returning to our case, if $M \setminus \Sigma$ turns out to be
subcritcal, the Floer homology $HF(\Gamma_L, \Gamma_L)$ vanishes.
Thus, starting with a Lagrangian submanifold $L \subset \Sigma$, we
have transformed it into a new Lagrangian $\Gamma_L \subset M
\setminus \Sigma$ whose Floer homology vanishes due to geometric
reasons. Note that in contrast to $\Gamma_L$, in general $L$ itself
cannot be Hamiltonianly displaced and in fact its Floer homology might
not vanish.

We now turn to algebraic computations in Floer homology. The idea is
to perform the computation of $HF(\Gamma_L, \Gamma_L)$ in an
alternative way. The main tool for this end is a spectral sequence
based on the theory developed by Oh~\cite{Oh:spectral,Oh:relative}.
The second step of this spectral sequence is the singular cohomology
of a Lagrangian, and the sequence converges to its Floer homology.
Comparing this computation (performed on $\Gamma_L$) with the
vanishing of $HF(\Gamma_L, \Gamma_L)$ we obtain restrictions on the
cohomology of $\Gamma_L$. In some cases we are even able to reproduce
the entire cohomology of $\Gamma_L$. Having done this, we derive
information on the cohomology of $L$ itself (recall that $\Gamma_L \to
L$ is a circle bundle). We give a rather detailed construction of this
spectral sequence in Section~\ref{S:Computations} since our approach
is somewhat different than Oh's original
work~\cite{Oh:spectral,Oh:relative}.  Detailed computations using the
spectral sequence appear throughout the proofs of the main theorems in
Section~\ref{S:Proofs-main}.

Let us turn now to the case when $M \setminus \Sigma$ is a critical
Stein manifold. In this case it is no longer true that all Lagrangian
submanifolds in $M \setminus \Sigma$ are displaceable. In
Section~\ref{Sb:Lag-trace} and~\ref{Sb:Displacable} we introduce a
kind of geometric obstruction for displaceability in $M \setminus
\Sigma$. This obstruction, which we call the {\em Lagrangian trace},
is a union of (possibly) immersed Lagrangian spheres $\Lambda \subset
\Sigma$. It has the property that for every Lagrangian $L \subset
\Sigma$ with $L \cap \Lambda = \emptyset$, the circle bundle $\Gamma_L
\subset M \setminus \Sigma$ can be Hamiltonianly displaced, in
particular $HF(\Gamma_L, \Gamma_L)=0$. Recomputing the vanishing of
$HF(\Gamma_L, \Gamma_L)$ using the spectral sequence we deduce that
Lagrangian submanifolds $L \subset \Sigma$ with certain topological
properties must intersect $\Lambda$.

\subsubsection*{The choice of the coefficients}
Our cohomological restrictions are with $\mathbb{Z}_2$-coefficients
only. This has to do with technical reasons coming from Floer homology
theory which is one of the main tools used in this paper.  However,
recent developments in Floer theory, due to Fukaya, Oh, Ohta and
Ono~\cite{FO3} make it possible to define, in some cases, Floer
homology with coefficients in $\mathbb{Q}$ or even $\mathbb{Z}$.  It
seems very likely that several of our results above continue to hold
for cohomology with $\mathbb{Z}$-coefficients as well.

\subsubsection*{The title of the paper}
``Lagrangian non-intersections'' is derived from the following idea
motivated in this paper: whenever the principle of Lagrangian
intersections fails (in the sense that a Lagrangian can be
Hamiltonianly displaced) we obtain restrictions on the topology of the
Lagrangian via computations in Floer homology.

\subsection*{Organization of the paper}
The rest of the paper is organized as follows.  In
Section~\ref{S:Stein} we collect some important facts from the
symplectic theory of Stein manifolds that will be used in the sequel.
We also develop in this section methods to displace Lagrangian
submanifolds in both subcritical and critical Stein manifolds.

In Sections~\ref{S:Polarizations} and~\ref{S:Lag-pol} we discuss
symplectic manifolds $\Sigma$ that appear as hyperplane sections in
other manifolds $M$. For Lagrangian submanifolds $L \subset \Sigma$ we
introduce the Lagrangian circle bundle construction giving rise to a
new Lagrangian submanifold $\Gamma_L \subset M \setminus \Sigma$. We
then study the possibilities to displace $\Gamma_L$ in $M \setminus
\Sigma$ and introduce the Lagrangian trace $\Lambda \subset \Sigma$
which is a kind of obstruction for displacing $\Gamma_L$. In
Section~\ref{Sb:trace_examples} we present explicit calculations of
$\Lambda$ for various examples of $\Sigma$.
Section~\ref{S:Computations} is devoted to computations in Floer
homology.  In Section~\ref{S:Proofs-main} we give the proofs of the
main theorems.  Finally, in Section~\ref{S:Generalizations} we present
some generalizations of the theorems of Section~\ref{S:Intro}.

\section{Symplectic geometry of Stein manifolds} \cntrs
\label{S:Stein}

Here we briefly recall some basic facts on Stein manifolds from the
symplectic viewpoint. The reader is referred to~\cite{El:Stein,El:psh}
for the foundations of the symplectic theory of Stein manifolds.
Apart from the contents of Subsection~\ref{Sb-Cocritical} most of the
material below can be found in~\cite{El-Gr:convex, El:Stein, El:psh,
  Bi-Ci:Stein}.

A {\em Stein manifold} is a triple $(V,J,\varphi)$ where $(V,J)$ is an
open complex manifold and $\varphi:V \to \mathbb{R}$ is a smooth
exhausting plurisubharmonic function. The term ``exhausting'' means
that $\varphi$ is proper and bounded from below. ``Plurisubharmonic''
means that the 2-form $\omega_{\varphi} = - dd^{\mathbb{C}}\varphi$ is
a $J$-positive symplectic form, i.e.~$-dd^{\mathbb{C}}\varphi(v,Jv)>0$
for every $0\neq v\in T(V)$. (Here and in what follows
$d^{\mathbb{C}}$ stands for operator that takes a smooth function
$\varphi$ to the $1$-form $d \varphi \circ J$.) We denote by
$g_{\varphi}(\cdot,\cdot)=\omega_{\varphi}(\cdot,J\cdot)$ the
associated Riemannian K\"{a}hler metric.

Given a Stein manifold $(V,J,\varphi)$ we have the gradient vector
field $X_{\varphi}=\textnormal{grad}_{g_{\varphi}}\varphi$ of
$\varphi$ with respect to the metric $g_{\varphi}$.  A simple
computation shows that $L_{X_\varphi}\omega_\varphi=\omega_\varphi$,
hence the flow $X_\varphi^t$ of $X_{\varphi}$ is conformally
symplectic, $(X_{\varphi}^t)^* \omega_{\varphi} = e^t
\omega_{\varphi}$. We remark that in contrast to other texts
(e.g.~\cite{El-Gr:convex}), we do not assume that the flow of
$X_{\varphi}$ is complete, unless explicitly stated. Note however that
since $\varphi$ is exhausting the flow $X_\varphi^t$ does exist for
all negative times.

\subsection{Canonical symplectic structures on Stein manifolds} \cntrsb
\label{Sb-Canonical}

Given a Stein manifold $(V, J, \varphi)$ and $R \in \mathbb{R}$, we
denote by $V_{\varphi\leq R}$ the sublevel set $\varphi^{-1}((-\infty,
R])$. We write $\textnormal{Crit}(\varphi)$ for the set of critical
points of the function $\varphi$.

Following~\cite{El-Gr:convex} we say that a Stein manifold
$(V,J,\varphi)$ is {\em complete} if the flow of gradient vector field
$X_{\varphi}=\textnormal{grad}_{g_{\varphi}}\varphi$ exists for all
positive times.

\begin{lem}[See~\cite{El-Gr:convex},~\cite{Bi-Ci:Stein}] 
   \label{L:completion}
   Let $(V,J,\varphi)$ be a Stein manifold. Then for every $R \in
   \mathbb{R}$ there exists an exhausting plurisubharmonic function
   $\varphi_R:V\to \mathbb{R}$ with the following properties:
   \begin{enumerate}
     \item $\varphi_R=\varphi$ on $V_{\varphi\leq R}$.
     \item $(V,J,\varphi_R)$ is a complete Stein manifold.
     \item $\textnormal{Crit}(\varphi_R) = \textnormal{Crit}(\varphi)$
      and for every $p \in \textnormal{Crit}(\varphi_R)$,
      $\textnormal{ind}_p(\varphi_R) = \textnormal{ind}_p(\varphi)$.
   \end{enumerate}   
   In particular, the inclusion $(V_{\varphi\leq R},\omega_{\varphi})
   \subset (V, \omega_{\varphi_R})$ is a symplectic embedding.
\end{lem}

The next lemma shows that the symplectic structure of a complete Stein
manifold is unique up to symplectomorphism.
\begin{lem}[See~\cite{El-Gr:convex}, compare~\cite{Bi-Ci:Stein}] 
   \label{L:canonical}
   Let $(V,J)$ be a complex manifold endowed with two exhausting
   plurisubharmonic functions $\varphi_1, \varphi_2$ such that both
   Stein manifolds $(V, J, \varphi_1)$ and $(V, J, \varphi_2)$ are
   complete.  Then the symplectic manifolds $(V, \omega_{\varphi_1})$
   and $(V, \omega_{\varphi_2})$ are symplectomorphic.
\end{lem}

\begin{rem} \label{R:completion}
   In view of this lemma, we shall sometimes denote the symplectic
   manifold associated to the completion of a Stein manifold $(V, J,
   \varphi)$ by $(V, \widehat{\omega})$ (since it does not depend on
   the choice of the plurisubharmonic function). Note that for every
   open subset $W \subset V$ with compact closure we have a symplectic
   embedding $(W, \omega_{\varphi}) \hookrightarrow (V,
   \widehat{\omega})$.
\end{rem}

\subsection{The skeleton of a Stein manifold $(V,J,\varphi)$} \cntrsb
\label{Sb-Skeleton}
By this we mean the subset $\Delta_{\varphi}\subset V$ which is formed
from the union of the stable submanifolds of the flow $X_{\varphi}^t$,
namely:
$$\Delta_{\varphi} = \bigcup_{p\in \textnormal{Crit}(\varphi)}
W^s_p(X_{\varphi}) = \left\{ x \in V \bigm| \lim_{t\to \infty}
   X_{\varphi}^t(x) \in \textnormal{Crit}(\varphi) \right\} .$$
\begin{lem} \label{L:cells}
   Let $(V,J,\varphi)$ be a Stein manifold and assume that $\varphi$
   is a Morse-Bott function. Then for every critical submanifold $C
   \subset V$ of $\varphi$ we have:
   \begin{enumerate}
     \item $C$ is isotropic with respect to $\omega_{\varphi}$.
     \item $\textnormal{ind}_C(\varphi) + \dim C \leq
      \frac{1}{2}\dim_{\mathbb{R}}V$.
     \item For every $p \in C$ the stable submanifold
      $W^s_p(X_{\varphi})$ is isotropic with respect to
      $\omega_{\varphi}$.
     \item For every $p\in C$ the unstable submanifold
      $W^u_p(X_{\varphi})$ is coisotropic with respect to
      $\omega_{\varphi}$.
   \end{enumerate}
\end{lem}
\begin{proof}
   For $\varphi$ being Morse a proof can be found in~\cite{El:psh}
   (see also~\cite{El-Gr:convex}, and see~\cite{Bi:barriers}
   Section~8). The case of Morse-Bott $\varphi$ is an obvious
   extension of the ``Morse case''.
\end{proof}

\begin{lem}[See~\cite{Bi:barriers}] \label{L:Morse}
   Let $(V,J,\varphi)$ be a Stein manifold, and assume that all the
   critical points of $\varphi$ lie in the subset $\{ \varphi < R \}$
   for some $R \in \mathbb{R}$. Then arbitrarily close to $\varphi'$
   in the $C^2$-topology there exists an exhausting plurisubharmonic
   function $\varphi':V \to \mathbb{R}$ such that:
   \begin{enumerate}
     \item $\varphi'=\varphi$ on $\{ \varphi \geq R \}$.
     \item $\varphi'$ is Morse.
     \item The flow of $X_{\varphi'}$ is Morse-Smale. In particular
      all trajectories of $X_{\varphi'}$ go from critical points of
      $\varphi'$ to either critical points of strictly higher index or
      to ``infinity'' (i.e. do not go to any other critical point).
      Moreover the skeleton $\Delta_{\varphi'}$ is an isotropic
      CW-complex (see below).
   \end{enumerate}
\end{lem}
Let $(Y, \omega)$ be a symplectic manifold and $\Delta \subset Y$ a
subset. We call $\Delta$ an {\em isotropic CW-complex} if there exists
an abstract CW-complex $K$ and a homeomorphism $i: K \to \Delta
\subset Y$ such that for every cell $C \subset K$ the restriction
$i|_{\textnormal{Int\,} C \approx \textnormal{Int\,} (D^{\dim C})} :
\textnormal{Int\,} C \to (Y,\omega)$ is an isotropic embedding.  We
refer the reader to~\cite{Bi:barriers} for more details on this
notion.

\subsection{Subcritical Stein manifolds} \cntrsb
\label{Sb-Subcritical}
Let $(V, J, \varphi)$ be a Stein manifold.  It is well
known~\cite{El:psh} that if $\varphi$ is Morse then for every critical
point $p$ we have $\textnormal{ind}_p \varphi \leq \frac{1}{2}
\dim_{\mathbb{R}} V$.  We say that $(V,J,\varphi)$ is {\em
  subcritical} if $\varphi$ is Morse with finite number of critical
points and for every $p \in \textnormal{Crit} \varphi$,
$\textnormal{ind}_p \varphi < \frac{1}{2} \dim_{\mathbb{R}} V$.  Note
that in this case $\dim \Delta_{\varphi} < \frac{1}{2}
\dim_{\mathbb{R}} V$ (hence the skeleton does not contain Lagrangian
cells).

The following lemma shows that in subcritical Stein manifolds any
compact subset can be Hamiltonianly displaced. See~\cite{Bi-Ci:Stein}
for the proof.
\begin{lem} \label{L:displacement-subcrit}   
   Let $(V,J,\varphi)$ be a complete subcritical Stein manifold.  Then
   for every compact subset $A\subset V$ there exists a compactly
   supported Hamiltonian diffeomorphism $h: (V,\omega_{\varphi}) \to
   (V,\omega_{\varphi})$ such that $h(A)\cap A = \emptyset$.
\end{lem}

\subsection{The critical coskeleton} \cntrsb
\label{Sb-Cocritical}
Let $(V, J, \varphi)$ be a Stein manifold, and suppose that $\varphi$
is a Morse-Bott function with finitely many critical submanifolds. We
denote by $p_1, \ldots, p_N$ the isolated critical points of $\varphi$
(if there are any). We say that $\varphi$ has property
$(\mathcal{S}_0)$ if one of the following two conditions is satisfied:
\begin{enumerate}
  \item For every positive dimensional critical submanifold $S$ of
   $\varphi$, $\textnormal{ind}_{S}(\varphi) + \dim S <
   \frac{1}{2}\dim_{\mathbb{R}} V$.
  \item $\varphi$ has no isolated critical points and only one
   positive dimensional critical manifold $S$, with $\dim S =
   \frac{1}{2} dim_{\mathbb{R}} V$.
\end{enumerate}
In the first case denote by $\{ p'_1, \ldots, p'_r \} \subset \{ p_1,
\ldots, p_N\}$ those critical points with
$\textnormal{ind}_{p'_i}(\varphi) = \frac{1}{2} \dim_{\mathbb{R}} V$
(again, it may happen that $r=0$.) In the second case pick a point
$p'_1 \in S$ and put $r=1$. We define the critical coskeleton
$\nabla_{\varphi}^{\textnormal{crit}} \subset V$ to be the union of
the unstable submanifolds of the $p'_i$'s, namely:
$$\nabla_{\varphi}^{\textnormal{crit}} = \bigcup_{i=1}^r W^u_{p'_i}
(X_{\varphi}) = \left\{ x \in V \bigm| \lim_{t\to -\infty}
   X_{\varphi}^t(x) \in \{ p'_1, \ldots, p'_r\} \right\} .$$

\begin{remnonum}
   Property $(\mathcal{S}_0)$ is purely technical and may look
   somewhat artificial. Note that if $\varphi$ is Morse then it
   automatically has property $(\mathcal{S}_0)$ (since all its
   critical points are isolated). Property $(\mathcal{S}_0)$ was
   created to accommodate a slightly more general situation than that.
   It covers two different (and unrelated) possibilities. The first
   possibility means that among the unstable submanifolds, those that
   are Lagrangian (i.e. have minimal dimension) all come from isolated
   critical points. The second possibility, roughly speaking, means
   that $(V,\omega_{\varphi})$ looks like a neighbourhood of the zero
   section in $T^*(S)$.
\end{remnonum}

\begin{exsnonum}
   \begin{enumerate}
     \item If $(V, J, \varphi)$ is a subcritical Stein manifold then
      clearly $r=0$, hence $\nabla_{\varphi}^{\textnormal{crit}} =
      \emptyset$.
     \item Let $M$ be a closed manifold and $V=T^*(M)$ be its
      cotangent bundle. Denote by $q \in M$ local coordinates along
      $X$ and by $p \in T_q^*(M)$ the dual coordinates along the
      cotangent fibres.  It is well known (see~\cite{El:psh}) that $V$
      can be endowed with the structure of a Stein manifold with
      $\varphi (q,p)= |p|^2$ and $X_{\varphi} =
      p\frac{\partial}{\partial p}$. Here $|\cdot|$ is a norm along
      the cotangent fibers (coming from a Riemannian metric on $M$).
      In this case the only critical submanifold is the zero section,
      and the critical coskeleton is just one fibre,
      $\nabla_{\varphi}^{\textnormal{crit}} = T_q^*(M)$.
   \end{enumerate}
\end{exsnonum}

\subsubsection{Property $(\mathcal{S})$} \label{sbsb:cond_S}
Let $(V,J,\varphi)$ be a Stein manifold. We say that $\varphi$ has
property $(\mathcal{S})$ if it has property $(\mathcal{S}_0)$ above
and in addition for every $1 \leq i \leq r$ all gradient trajectories
of $X_{\varphi}$ emanating from the points $p'_i$ go to ``infinity''.
Note that this definition does not depend on the choice of the point
$p'_1$ in the case when $\varphi$ has only one critical submanifold
$S$ with $\dim S = \frac{1}{2}\dim_{\mathbb{R}}V$. Indeed, in that
case $S$ is the minimum of $\varphi$, hence all gradient trajectories
emanating from points of $S$ go to ``infinity''.

\begin{lem} \label{L:displacement-crit}
   Let $(V,J,\varphi)$ be a complete Stein manifold and assume that
   $\varphi$ has property $(\mathcal{S})$. Let $A \subset V$ be a
   compact subset with $A \cap \nabla_{\varphi}^{\textnormal{crit}} =
   \emptyset$. Then there exists a compactly supported Hamiltonian
   diffeomorphism $h:(V, \omega_{\varphi}) \to (V, \omega_{\varphi})$
   such that $h(A) \cap A = \emptyset$.
\end{lem}

\begin{exsnonum}
   \begin{enumerate}
     \item If $(V, J, \varphi)$ is subcritical then the the lemma
      reduces to Lemma~\ref{L:displacement-subcrit} since in this case
      $\nabla_{\varphi}^{\textnormal{crit}} = \emptyset$, hence any
      compact subset can be Hamiltonianly displaced.
     \item Let $V=T^*(M)$. As we have just seen above
      $\nabla_{\varphi}^{\textnormal{crit}} = T^*_q(M)$ for some $q
      \in M$. Hence we recover a statement due to Lalonde and
      Sikorav~\cite{La-Si:Lag} that any compact subset of $T^*(M)$
      lying in the complement of a fibre $T^*_q(M)$ can be
      Hamiltonianly displaced.
   \end{enumerate}
\end{exsnonum}

Before we prove Lemma~\ref{L:displacement-crit} we shall need some
preparations. Given a Morse function $\varphi : V \to \mathbb{R}$,
denote by $\textnormal{Crit}_{\leq k}(\varphi)$ the set of critical
points of $\varphi$ of index $\leq k$. Denote by $\Delta_{\varphi}^k$
the subskeleton
$$\Delta_{\varphi}^k = \bigcup_{p\in \textnormal{Crit}_{\leq
    k}(\varphi)} W^s_p(X_{\varphi}) = \left\{ x \in V \bigm|
   \lim_{t\to \infty} X_{\varphi}^t(x) \in \textnormal{Crit}(\varphi)
\right\}.$$

We shall need the following Proposition for the proof of
Lemma~\ref{L:displacement-crit}.

\begin{prop} \label{P:Delta-k}
   Let $(V,J,\varphi)$ be a Stein manifold. Fix an integer $0\leq k
   \leq \dim_{\mathbb{C}}V$. Assume that:
   \begin{enumerate}
     \item $\varphi$ is a Morse function with finitely many critical
      points $x_1, \ldots, x_{\nu} \in V$, arranged so that
      $\textnormal{Crit}_{\leq k}(\varphi)=\{x_1, \ldots, x_l\}$,
      $l\leq \nu$.
     \item There are no trajectories of $X_{\varphi}$ that connect any
      of the critical points $x_{l+1}, \ldots, x_{\nu}$ with one of
      the critical points $x_1, \ldots, x_l$.
   \end{enumerate}
   Fix mutually disjoint neighbourhood $U_1, \ldots, U_l$ of $x_1,
   \ldots, x_l$ respectively. Then arbitrarily close to $\varphi$ in
   the $C^2$-topology there exists an exhausting plurisubharmonic
   function $\varphi'$ with the following properties:
   \begin{enumerate}
     \item $\varphi=\varphi'$ on $V \setminus (U_1 \cup \ldots \cup
      U_l)$.
     \item $\varphi'$ is Morse,
      $\textnormal{Crit}(\varphi')=\textnormal{Crit}(\varphi)$ and for
      every $1\leq i \leq \nu$,
      $\textnormal{ind}_{x_i}(\varphi')=\textnormal{ind}_{x_i}(\varphi)$.
     \item The flow of $X_{\varphi'}$ connects any of the critical
      points $x_1, \ldots, x_l$ either to a point of strictly higher
      index or to ``infinity''. Moreover the subskeleton
      $\Delta_{\varphi'}^{k}$ is an isotropic CW-complex.
   \end{enumerate}
\end{prop}
We omit the proof as it is a straightforward adaptation of the
arguments from Section~9 of~\cite{Bi:barriers}.

\begin{proof}[Proof of Lemma~\ref{L:displacement-crit}]
   The proof generalizes ideas from~\cite{Bi-Ci:Stein} (see Lemma~3.2
   there).
   
   Assume first that $\varphi$ is a Morse function, hence it has no
   positive dimensional critical submanifolds. Put $k=\frac{1}{2}
   \dim_{\mathbb{R}} V - 1$. Applying Proposition~\ref{P:Delta-k} we
   obtain a new plurisubharmonic function $\varphi'$ for which
   $\Delta_{\varphi'}^{k}$ is a CW-complex. Note that if we choose the
   neighborhoods $U_1, \ldots, U_l$ of the points in
   $\textnormal{Crit}_{\leq k}(\varphi)$ to be small enough we can
   arrange that $\nabla_{\varphi'}^{\textnormal{crit}} =
   \nabla_{\varphi}^{\textnormal{crit}}$.
   
   By Moser argument there is a symplectomorphism $f:(V,
   \omega_{\varphi}) \to (V, \omega_{\varphi'})$ which is supported in
   $U_1 \cup \ldots \cup U_l$. In particular
   $f(\nabla_{\varphi}^{\textnormal{crit}}) =
   \nabla_{\varphi'}^{\textnormal{crit}}$, and $f(A) \cap
   \nabla_{\varphi'}^{\textnormal{crit}} = \emptyset$. Thus by
   replacing $\varphi$ by $\varphi'$ we may assume without loss of
   generality that $\Delta_{\varphi}^{k}$ is a CW-complex.
   
   Since $\dim \Delta_{\varphi}^{k} < \frac{1}{2} \dim_{\mathbb{R}} V$
   there exists a Hamiltonian isotopy $g_t:(V,\omega_{\varphi}) \to
   (V,\omega_{\varphi})$, compactly supported in an arbitrarily small
   neighbourhood of $\Delta_{\varphi}^k$, such that
   $g_1(\Delta_{\varphi}^k) \cap \Delta_{\varphi}^k = \emptyset$. As
   $\Delta_{\varphi}^k$ is compact there exists a small neighbourhood
   $W$ of $\Delta_{\varphi}^k$ so that $g_1(W)\cap W = \emptyset$.
   
   Since $A \cap \nabla_{\varphi}^{\textnormal{crit}} = \emptyset$,
   for large enough $T>0$ we have $X_{\varphi}^{-T}(A) \subset W$.  As
   $g_1$ displaces $W$ we have:
   $$X_{\varphi}^{T} \circ g_1 \circ X_{\varphi}^{-T}(A) \cap A =
   \emptyset.$$
   Finally, it is a straightforward computation to check
   that
   $$h_t = X_{\varphi}^{tT}\circ g_t \circ X_{\varphi}^{-tT}$$
   is a
   Hamiltonian isotopy (see~\cite{Bi-Ci:Stein}, Lemma~3.2).

   Assume now that $\varphi$ has also positive dimensional
   submanifolds say $S_1, \ldots, S_q$ with
   $\textnormal{ind}_{S_i}(\varphi) + \dim S_i < \frac{1}{2}
   \dim_{\mathbb{R}}V$ for every $i$. Pick for every $i$ a generic
   Morse function $f_i: S_i \to \mathbb{R}$ and a cut off function
   $\rho_i$ which is identically $1$ near $S_i$ and identically $0$
   outside a small neighbourhood $W_i$ of $S_i$. Consider now the
   function $\varphi_{\epsilon} = \varphi + \epsilon \sum_{i=1}^q
   \rho_i f_i$. Clearly for small $\epsilon$, $\varphi_{\epsilon}$ is
   plurisubharmonic.  Moreover $\varphi_{\epsilon}$ is Morse and its
   critical points consists of the isolated critical points of
   $\varphi$ and the critical points of the Morse functions $f_1,
   \ldots, f_q$.  Moreover, for every critical point $p \in
   \textnormal{Crit}(f_i)$, we have $\textnormal{ind}_p
   (\varphi_{\epsilon}) = \textnormal{ind}_p(f_i) +
   \textnormal{ind}_{S_i}(\varphi)$. From property $(\mathcal{S}_0)$
   we get $\textnormal{ind}_p (\varphi_{\epsilon}) < \frac{1}{2}
   dim_{\mathbb{R}} V$.  Therefore the critical points of
   $\varphi_{\epsilon}$ of index $\frac{1}{2}\dim_{\mathbb{R}} V$ are
   exactly the same as those of $\varphi$. Moreover
   $\varphi_{\epsilon} = \varphi$ near these points. Next, note that
   if the perturbation above is in small enough neighborhoods $W_i$ of
   the $S_i$'s then due to assumption $(\mathcal{S})$ we have
   $\nabla_{\varphi_{\epsilon}}^{\textnormal{crit}} =
   \nabla_{\varphi}^{\textnormal{crit}}$.
   
   Now, by Moser argument there is a symplectomorphism $f:(V,
   \omega_{\varphi}) \to (V, \omega_{\varphi_{\epsilon}})$ which is
   supported in $W_1 \cup \ldots \cup W_q$. And again,
   $f(\nabla_{\varphi}^{\textnormal{crit}}) =
   \nabla_{\varphi_{\epsilon}}^{\textnormal{crit}}$, and $f(A) \cap
   \nabla_{\varphi_{\epsilon}}^{\textnormal{crit}} = \emptyset$.
   Replacing $\omega_{\varphi}$ by $\omega_{\varphi_{\epsilon}}$ and
   $A$ by $f(A)$ we arrive to the case from the beginning of the
   proof.
   
   Finally assume that the only critical points of $\varphi$ consist
   of one critical submanifolds $S$ of dimension $\frac{1}{2}
   \dim_{\mathbb{R}} V$. Consider the map $F: V \to S$ defined by
   $$F(x) = \lim_{t \to -\infty} X_{\varphi}^t (x).$$
   Since $\varphi$
   is Morse-Bott this map is a locally trivial fibration in a
   neighbourhood of $S$ (see~\cite{Au-Br}). Note that
   $\nabla_{\varphi}^{\textnormal{crit}}$ is just the preimage under
   $F$ of a point $p \in S$. Therefore since $A \cap
   \nabla_{\varphi}^{\textnormal{crit}} = \emptyset$ we have $F(A)
   \subset S \setminus \{ p \}$.
   
   Since $S$ is Lagrangian a small neighbourhood of $S$ can be
   identified with a neighbourhood of $T^*(S)$. Pick a Morse function
   $G$ on $S$ having all its critical points in $S \setminus F(A)$.
   Clearly the Hamiltonian isotopy $g_t$ of $G$ (viewed as a
   Hamiltonian in $T^*(S)$) will displace $F(A)$ away of $S$ within
   arbitrary small time, say $t=\epsilon$ (compare~\cite{La-Si:Lag}).
   The result now follows in the same way as in the beginning of the
   proof. Namely for large enough $T>0$ the Hamiltonian diffeomorphism
   $$X_{\varphi}^{T} \circ g_{\epsilon} \circ X_{\varphi}^{-T}(A) \cap
   A = \emptyset$$
   will displace $A$.
\end{proof}

\section{Polarizations and decompositions of K\"{a}hler manifolds} \cntrs
\label{S:Polarizations}
A basic tool that we shall use throughout this work is a decomposition
technique for K\"{a}hler manifolds that was developed
in~\cite{Bi:barriers}.  In
Subsections~\ref{Sb:pol}-~\ref{Sb:subcritical_polarizations} we
briefly summarize the necessary facts from~\cite{Bi:barriers} where
more details can be found. In Subsection~\ref{Sb:Lag-trace} we
introduce the concept of Lagrangian trace and
in~\ref{Sb:trace_examples} we compute some examples.

\subsection{Polarized K\"{a}hler manifolds} \cntrsb \label{Sb:pol}
Throughout this paper, by a K\"{a}hler manifold we mean a triple $(M,
\omega, J)$ where $(M, \omega)$ is a closed symplectic manifold and
$J$ is an (integrable) complex structure compatible with $\omega$.

A {\em polarized K\"{a}hler manifold} $\mathcal{P} = (M^{2n}, \omega,
J; \Sigma)$ is a K\"{a}hler manifold $(M, \omega, J)$ with $[\omega]
\in H^2(M;\mathbb{Z})$ together with a smooth and reduced complex
hypersurface $\Sigma \subset M$ whose homology class $[\Sigma] \in
H_{2n-2}(M)$ is the Poincar\'{e} dual to $k[\omega] \in H^2(M)$ for
some $k \in \mathbb{N}$. The number $k$ will be called {\em the degree
  of the polarization} $\mathcal{P}$ and denoted by $k_{\mathcal{P}}$.
Note that our notion of polarized K\"{a}hler manifolds is slightly
different from the one common in algebraic geometry.

\subsection{Additional structures associated with a polarization} 
\cntrsb \label{Sb:additinal}

We shall now define a distinguished plurisubharmonic function
$\varphi_{_{\mathcal{P}}}: M\setminus\Sigma \to \mathbb{R}$ which is
canonically associated with the polarization $\mathcal{P}$.  For this
purpose let $\mathcal{L} = \mathcal{O}_M(\Sigma)$ be the holomorphic
line bundle defined by the divisor $\Sigma$. Denote by $s: M \to
\mathcal{L}$ the (unique up to a constant factor) holomorphic section
whose zero set $\{ s=0 \}$ is $\Sigma$. Choose a hermitian metric $\|
\cdot \|$ on $\mathcal{L}$, and a compatible connection $\nabla$ with
curvature $R^{\nabla} = 2\pi i k_{\mathcal{P}} \omega$.  Finally,
define $\varphi_{_{\mathcal{P}}} : M\setminus\Sigma \to \mathbb{R}$ to
be
$$\varphi_{_{\mathcal{P}}}(x) = -\frac{1}{4\pi k_{\mathcal{P}}} \log
\|s(x)\|^2.$$
Put $V=M\setminus \Sigma$.  A simple computation shows
that $-dd^{\mathbb{C}}\varphi_{_{\mathcal{P}}} = \omega$, hence
$\varphi_{_{\mathcal{P}}}$ is plurisubharmonic.  Moreover, it is not
hard to see that $\varphi_{_{\mathcal{P}}}$ is exhausting and that it
has no critical points outside some compact subset of $V$
(see~\cite{Bi:barriers}).

It is important to remark that the function $\varphi_{_{\mathcal{P}}}$
is canonically determined by the polarization $\mathcal{P}$ up to an
additive constant and does not depend on any of the choices made for
$\|\cdot \| , \, s$ or $\nabla$. This is due to the requirement on the
curvature $R^{\nabla}$ and the fact that $J$ is integrable
(see~\cite{Bi:barriers} for more details).  Next, let $g_{_{\omega,
    J}} = \omega(\cdot, J \cdot)$ be the K\"{a}hler Riemannian metric
associated with the pair $(\omega,J)$. Finally denote by
$X_{_{\mathcal{P}}}^t$ the gradient flow of $\varphi_{_{\mathcal{P}}}$
with respect to $g_{_{\omega,J}}$. (Note that $X_{_{\mathcal{P}}}^t$
is not complete for $t>0$, since $(V,\omega)$ has finite volume.)

Consider the Stein manifold $(V = M \setminus \Sigma, J,
\varphi_{_{\mathcal{P}}})$. We denote by $\Delta_{\mathcal{P}} \subset
V$ its skeleton (see~\ref{Sb-Skeleton} above).  Note that
$\Delta_{\mathcal{P}}\subset M\setminus \Sigma$ is compact since the
flow $X_{_{\mathcal{P}}}^t$ is complete at $-\infty$ and
$\textnormal{Crit}(\varphi_{_{\mathcal{P}}})$ is a compact subset of
$M\setminus \Sigma$. We remark that $\Delta_{\mathcal{P}}$ is
completely determined by the polarization $\mathcal{P}$ without any
further choices since the function $\varphi_{_{\mathcal{P}}}$ is
determined (up to an additive constant) by $\mathcal{P}$.  We shall
therefore call $\Delta_{\mathcal{P}}$ {\em the skeleton} associated
with the polarization $\mathcal{P}$.

\subsection{The decomposition associated to a polarization} \cntrsb
\label{Sb:decompostion}
In this section we explain how to decompose a K\"{a}hler manifold into
two basic building blocks. The first piece is a standard symplectic
disc bundle over a complex hypersurface. The second piece is the
isotropic skeleton of the Stein manifold which is the complement of
this hypersurface.

\subsubsection{Standard symplectic disc bundles} 
\label{sbsb:sdb}

Let $\mathcal{P} = (M, \omega, J; \Sigma)$ be a polarization of degree
$k_{\mathcal{P}}$ of a K\"{a}hler manifold.

Put $\omega_{_{\Sigma}}=\omega|_{T(\Sigma)}$ and let
$\pi:N_{\Sigma}\to \Sigma$ be the (complex) normal line bundle of
$\Sigma$ in $M$ with first Chern class $c_1^{N_{\Sigma}} =
k_{\mathcal{P}}[\omega_{_{\Sigma}}]\in H^2(\Sigma)$. Let $\|\cdot\|$
be any hermitian metric on $N_{\Sigma}$ and denote by $E_{\Sigma} = \{
v \in N_{\Sigma} \;\bigl|\; \|v\| < 1 \}$ the open unit disc bundle of
$N_{\Sigma}$.  Choose a connection $\nabla$ on $N_{\Sigma}$ with
curvature $R^{\nabla} = 2\pi i k_{\mathcal{P}}\omega_{_{\Sigma}}$ and
denote by $\alpha^{\nabla}$ the associated {\em transgression 1-form}
on $N_{\Sigma} \setminus 0$ defined by:
\begin{itemize}
  \item $\alpha^{\nabla}_{(u)}(u) = 0, \; \alpha^{\nabla}_{(u)}(iu) =
   \frac{1}{2\pi}$ for every $u \in N_{\Sigma} \setminus 0$.
  \item $\alpha^{\nabla}|_{H^{\nabla}} = 0$, where $H^{\nabla}$ is the
   horizontal distribution of $\nabla$.
\end{itemize}
With this normalization of $\alpha^{\nabla}$ we have $d\alpha^{\nabla}
= -\pi^* (k_{\mathcal{P}}\omega_{_{\Sigma}})$.  Define now the
following symplectic form $\omega_{\textnormal{can}}$ on $E_{\Sigma}$:
$$\omega_{\textnormal{can}} = k_{\mathcal{P}}\pi^*\omega_{_{\Sigma}} +
d(r^2 \alpha^{\nabla}),$$
where $r$ is the radial coordinate along the
fibres induced by $\| \cdot \|$. It is easy to check that
$\omega_{\textnormal{can}}$ is well defined, that it is symplectic,
and has the following three properties:
\begin{enumerate}
  \item All fibres of $\pi : E_{\Sigma} \to \Sigma$ are symplectic
   with respect to $\omega_{\textnormal{can}}$ and have area $1$.
  \item The restriction of $\omega_{\textnormal{can}}$ to the zero
   section $\Sigma \subset E_{\Sigma}$ equals
   $k_{\mathcal{P}}\omega_{_{\Sigma}}$.
  \item $\omega_{\textnormal{can}}$ is $S^1$-invariant with respect to
   the obvious circle action on $E_{\Sigma}$.
\end{enumerate}

Although $\omega_{\textnormal{can}}$ a priori depends on $\| \cdot \|$
and $\nabla$, different choices of these structures in fact lead to
symplectically equivalent results (see~\cite{McD-Sa:Intro}, see
also~\cite{Bi:barriers, Bi-Ci:closed}). We shall henceforth call
$(E_{\Sigma}, \omega_{\text{can}})$ {\em the standard symplectic disc
  bundle over $(\Sigma, \omega_{_{\Sigma}})$ modeled on $N_{\Sigma}$}.
Often we shall multiply $\omega_{\text{can}}$ by a positive number
$c>0$ (usually by $c=\frac{1}{k_{\mathcal{P}}}$) and refer to
$(E_{\Sigma}, c \omega_{\text{can}})$ as {\em the standard symplectic
  disc bundle with fibres of area $c$}. (Note that now the restriction
of this symplectic form to $\Sigma \subset E_{\Sigma}$ equals $c
k_{\mathcal{P}}\omega_{_{\Sigma}}$, not
$k_{\mathcal{P}}\omega_{_{\Sigma}}$.)

\subsubsection{The decomposition}
\label{sbsb:decompostion}
Let $\mathcal{P} = (M, \omega, J; \Sigma)$ be a polarized K\"{a}hler
manifold. Denote by $\rho_{_{\mathcal{P}}}: M \to \mathbb{R}$ the
function $\rho_{_{\mathcal{P}}}(x) = \|s(x)\|^2$, and let
$Z_{_{\mathcal{P}}}$ be the gradient vector field of
$-\rho_{_{\mathcal{P}}}$. Note that since $\rho_{_{\mathcal{P}}} =
e^{-4\pi k_{\mathcal{P}} \varphi_{_{\mathcal{P}}}}$ the vector fields
$Z_{_{\mathcal{P}}}$ and $X_{_{\mathcal{P}}}$ are positively
proportional on $M \setminus \Sigma$.

\begin{thm}[See~\cite{Bi:barriers}] \label{T:decomposition}
   Let $\mathcal{P} = (M, \omega, J; \Sigma)$ be a polarized
   K\"{a}hler manifold. Then, the complement of the skeleton $(M
   \setminus \Delta_{\mathcal{P}}, \omega)$ is symplectomorphic to the
   following standard symplectic disc bundle over $\Sigma$
   $$(E_{\Sigma}, \frac{1}{k_{\mathcal{P}}}\omega_{\text{can}}) \to
   (\Sigma, \, k_{\mathcal{P}} \omega_{_{\Sigma}})$$
   which is modeled
   on the normal bundle $N_{\Sigma}$, and with fibres of area
   $\frac{1}{k_{\mathcal{P}}}$. In fact, there exists a canonical
   symplectomorphism $F_{\mathcal{P}}$, which depends only on
   $\mathcal{P}$, such that the following diagram commutes:
   \[
      \begin{CD}
         (E_{\Sigma}, \frac{1}{k_{\mathcal{P}}}\omega_{\text{can}})
         @>{F_{\mathcal{P}}} >> (M \setminus \Delta_{\mathcal{P}},\omega) \\
         @A{\text{$0$-section}}AA @AA{\text{inclusion}}A \\
         (\Sigma, \omega_{_{\Sigma}}) @= (\Sigma, \omega_{_{\Sigma}})
      \end{CD}
   \]
   Moreover, $F_{\mathcal{P}}$ sends the flow lines of
   $Z_{_{\mathcal{P}}}$ to the lines of the negative radial flow on
   $E_{\Sigma}$, namely $D F_{\mathcal{P}}(Z_{_{\mathcal{P}}})$ is
   negatively proportional to the radial vector field
   $r\frac{\partial}{\partial r}$ on $E_{\Sigma}$.
\end{thm}
The proof of this theorem appears in~\cite{Bi:barriers}. The
``Moreover'' statement, is not stated explicitly in~\cite{Bi:barriers}
as a theorem but is proved there (see proof of Proposition~7.B
in~\cite{Bi:barriers}).

\subsection{Subcritical polarizations} \cntrsb
\label{Sb:subcritical_polarizations}
A polarization $\mathcal{P} = (M, \omega, J; \Sigma)$ is called {\em
  subcritical} if there exists a plurisubharmonic function
$\varphi:(V=M\setminus \Sigma, J) \to \mathbb{R}$ such that
$(V,J,\varphi)$ is a subcritical Stein manifold (namely $\varphi$ is
Morse and for every $p \in \textnormal{Crit}(\varphi)$,
$\textnormal{ind}_p(\varphi) < \frac{1}{2}\dim_{\mathbb{R}}V$).  Note
that we do not assume $\varphi$ to be the canonical function
$\varphi_{_{\mathcal{P}}}$. We refer the reader to~\cite{Bi-Ci:closed}
for more information on subcritical polarizations, examples and
criteria for identifying them.

\subsection{The Lagrangian trace} \cntrsb
\label{Sb:Lag-trace}
Let $\mathcal{P} = (M, \omega, J; \Sigma)$ be a polarized K\"{a}hler
manifold. We say that the polarization $\mathcal{P}$ has property
$(\mathcal{S})$ if $\varphi_{_{\mathcal{P}}}: M \setminus \Sigma \to
\mathbb{R}$ has property $(\mathcal{S})$ of
Section~~\ref{sbsb:cond_S}. We denote by
$\nabla_{\mathcal{P}}^{\textnormal{crit}} \subset M \setminus \Sigma$
the critical coskeleton of $\varphi_{_{\mathcal{P}}}$.

Let $\rho_{_{\mathcal{P}}}: M \to \mathbb{R}$ be the function defined
by $\rho_{_{\mathcal{P}}}(x) = \|s(x)\|^2$.  Note that the critical
points of $\rho_{_{\mathcal{P}}}$ consist of those of
$\varphi_{_{\mathcal{P}}}$ and $\Sigma$ which is a non-degenerate
critical submanifold (of index $0$). On $M \setminus \Sigma$ the
gradients of $\varphi_{_{\mathcal{P}}}$ and of
$-\rho_{_{\mathcal{P}}}$ have the same (oriented) flow lines. Denote
by $Z_{_{\mathcal{P}}}^t:M \to M$ the gradient flow of
$-\rho_{_{\mathcal{P}}}$.

Assume now that $\mathcal{P}$ has property $(\mathcal{S})$ and let
$\nabla_{\mathcal{P}}^{\textnormal{crit}} \subset M \setminus \Sigma$
be its critical coskeleton.  Denote by $\Lambda_{\mathcal{P}} \subset
\Sigma$ the subset obtained from
$\nabla_{\mathcal{P}}^{\textnormal{crit}}$ by ``projecting'' it using
$\lim_{t \to \infty} Z_{_{\mathcal{P}}}^t$ to $\Sigma$, namely
$$\Lambda_{\mathcal{P}} = \left\{ x \in \Sigma \bigm| x = \lim_{t \to
     \infty} Z_{_{\mathcal{P}}}^t(p), \textnormal{for some } p \in
   \nabla_{\mathcal{P}}^{\textnormal{crit}} \right\}.$$
In case
$\nabla_{\mathcal{P}}^{\textnormal{crit}} = \emptyset$ (e.g. if
$\varphi_{_{\mathcal{P}}}$ is subcritical) we put
$\Lambda_{\mathcal{P}} = \emptyset$.  We call
$\Lambda_{\mathcal{P}}\subset \Sigma$ the {\em Lagrangian trace} of
the polarization $\mathcal{P}$. This term is justified by the
following proposition.

\begin{prop} \label{P:Ltrace}
   Let $\mathcal{P} = (M, \omega, J; \Sigma)$ be a polarized
   K\"{a}hler manifold with $\varphi_{_{\mathcal{P}}}$ having property
   $(\mathcal{S})$. Let $p'_1, \ldots, p'_r$ be those critical points
   of $\varphi_{_{\mathcal{P}}}$ as chosen in~\ref{Sb-Cocritical}
   above. Then the corresponding Lagrangian trace
   $\Lambda_{\mathcal{P}} \subset \Sigma$ consists of a union of $r$
   immersed (but possibly embedded) Lagrangian spheres one for each of
   the points $p'_1, \ldots, p'_r$.
\end{prop}

\begin{proof}
   Denote by $G: M\setminus \Delta_{\mathcal{P}} \to \Sigma$ the end
   point map of the flow $Z_{_{\mathcal{P}}}^t$, namely $G(x) =
   \lim_{t \to \infty} Z_{_{\mathcal{P}}}^t(x)$. Clearly
   $$\Lambda_{\mathcal{P}} = \cup_{i=1}^r G(W'_i),$$
   where $W'_i =
   W_{p'_i}^u(X_{_{\mathcal{P}}}) \setminus \{ p'_i \} =
   W_{p'_i}^u(Z_{_{\mathcal{P}}}) \setminus \{ p'_i \}$ are the
   unstable submanifolds of $X_{_{\mathcal{P}}}$ at $p'_i$. We shall
   prove that for every $i$, $G(W'_i)$ is an immersed Lagrangian
   sphere in $\Sigma$.
   
   Denote by $\pi: E_{\Sigma} \to \Sigma$ the standard symplectic disc
   bundle, endowed with the symplectic structure
   $\frac{1}{k_{\mathcal{P}}} \omega_{\textnormal{can}}$. Using
   Theorem~\ref{T:decomposition} we may identify $(E_{\Sigma},
   \frac{1}{k_{\mathcal{P}}} \omega_{\textnormal{can}})$ with $(M
   \setminus \Delta_{\mathcal{P}}, \omega)$. We shall view from now on
   $W'_i$ as a submanifold of $E_{\Sigma}$ and $G$ as a map
   $G:E_{\Sigma} \to \Sigma$. It follows from
   Theorem~\ref{T:decomposition} that the map $G$ coincides with
   projection map $\pi: E_{\Sigma} \to \Sigma$.
   
   Denote $P_{\epsilon} = \{ v \in E_{\Sigma} \mid \| v \| = \epsilon
   \}$. Note that $W'_i$ intersects $P_{\epsilon}$ transversely
   because the vector field $\frac{\partial}{\partial r}$ is tangent
   to $W'_i$. Put $L_i = W'_i \cap P_{\epsilon}$. We claim that $L_i$
   is diffeomorphic to a sphere. Indeed, pick a small ball $B_i
   \subset W_{p'_i}^u$ centered around $p'_i$ whose boundary $\partial
   B_i$ is transverse to $Z_{_{\mathcal{P}}}$.  Taking $B_i$ to be
   small enough we may assume that every flow line of
   $Z_{_{\mathcal{P}}}$ intersects $\partial B_i$ exactly once.
   Denote by $\nu: E_{\Sigma} \setminus \Sigma \to P_{\epsilon}$ the
   map $\nu(v) = \epsilon \frac{v}{\| v \|}$. By our choice of $B_i$
   we have that $\nu$ sends the sphere $\partial B_i$
   diffeomorphically onto $L_i$.
   
   Next, we claim that $\pi|_{L_i}:L_i \to \Sigma$ is an immersion. To
   prove this, note that for every $x \in E_{\Sigma} \setminus
   \Sigma$, $\ker D\pi_x = \mathbb{R} \frac{\partial}{\partial r}
   \oplus i \mathbb{R} \frac{\partial}{\partial r}$. Now let $v = a
   \frac{\partial}{\partial r} + ib \frac{\partial}{\partial r} \in
   \ker D\pi_{x} \cap T_x(L_i)$. Since $L_i \subset P_{\epsilon}$, we
   have $a=0$. By Lemma~\ref{L:cells} $W'_i$ is Lagrangian (it is
   coisotropic and has half the dimension of $M$). As
   $\frac{\partial}{\partial r}$ is tangent to $W'_i$, we have
   $\omega_{\textnormal{can}}(\frac{\partial}{\partial r}, v) = 0$.
   But $\omega_{\textnormal{can}}(\frac{\partial}{\partial r}, v) = b
   \omega_{\textnormal{can}} (\frac{\partial}{\partial r}, i
   \frac{\partial}{\partial r})$ which can vanish only if $b=0$. Thus
   $v=0$. This proves that $\pi|_{L_i}$ is an immersion.
   
   It remains to prove that $\pi(L_i)$ is Lagrangian. Let $\xi = \ker
   (\alpha^{\nabla}|_{T(P_{\epsilon})})$ be the contact distribution
   on $P_{\epsilon}$. Then $T(P_{\epsilon}) = \xi \oplus i \mathbb{R}
   \frac{\partial}{\partial r}$. We first claim that $T(L_i) \subset
   \xi$. Indeed, let $v=u+i a \frac{\partial}{\partial r} \in T(L_i)$,
   where $u \in \xi$, $a \in \mathbb{R}$. As before
   $\omega_{\textnormal{can}}(v, \frac{\partial}{\partial r}) = 0$,
   because $W'_i$ is Lagrangian. But $\omega_{\textnormal{can}}(v,
   \frac{\partial}{\partial r}) = a \omega_{\textnormal{can}}(i
   \frac{\partial}{\partial r}, \frac{\partial}{\partial r})$ which
   can vanish only if $a=0$. Thus $v=u \in \xi$. Finally note that
   $$\omega_{\textnormal{can}} = k_{\mathcal{P}} \pi^*
   \omega_{_{\Sigma}} + d(r^2 \alpha^{\nabla}) =
   k_{\mathcal{P}}(1-r^2)\omega_{_{\Sigma}} + 2r dr \wedge
   \alpha^{\nabla},$$
   hence $D\pi$ sends $(\xi,
   (1-\epsilon^2)\omega_{\textnormal{can}}|_{\xi})$ isomorphically to
   $(T(\Sigma), \omega_{_{\Sigma}})$. As $T(L_i)$ is Lagrangian in
   $(\xi, \omega_{\textnormal{can}}|_{\xi})$ we conclude that
   $D\pi(T(L_i))$ is Lagrangian in $(T(\Sigma), \omega_{_{\Sigma}})$.
\end{proof}

\subsection{Examples} \cntrsb \label{Sb:trace_examples} 
Let us present a few explicit examples of Lagrangian traces coming
from various polarizations.

\subsubsection{Subcritical polarizations} \label{sbsb:subcritical}
Let $\mathcal{P} = (M, \omega, J; \Sigma)$ be a subcritical
polarization (see~\ref{Sb:subcritical_polarizations}). In this case
$\nabla_{\varphi}^{\textnormal{crit}} = \emptyset$, hence
$\Lambda_{\mathcal{P}} = \emptyset$.

The simplest example of a subcritical polarization consists of
$M={\mathbb{C}}P^n$ endowed with its standard symplectic K\"{a}hler
form $\sigma$ and $\Sigma \approx {\mathbb{C}}P^{n-1}$ being a linear
hyperplane. The skeleton in this case is a point $\Delta_{\mathcal{P}}
= \textnormal{pt}$ (see~\cite{Bi:barriers} for more details).

Another example of a subcritical polarization is $M = {\mathbb{C}}P^n
\times {\mathbb{C}}P^{n+r}$, $r \geq 1$, endowed with the split
symplectic structure $\sigma \oplus \sigma$, and
$$\Sigma = \left\{ ([z_0: \ldots : z_n], [w_0: \ldots :w_{n+r}]) \in
   {\mathbb{C}}P^n \times {\mathbb{C}}P^{n+r} \Bigm| \sum_{i=0}^{n-1}
   z_i w_i = z_n \sum_{j=n}^{n+r}w_j \right\}.$$
In this case the
skeleton turns out to be an isotropic copy of ${\mathbb{C}}P^n$
(see~\cite{Bi:barriers} for more details).

We refer the reader to~\cite{Bi-Ci:closed} for more details and
examples of subcritical polarizations.

\subsubsection{The quadric} \label{sbsb:quadric}
Consider the polarization $\mathcal{P} = (M, \omega, J; \Sigma)$ with
$M={\mathbb{C}}P^{n+1}$ and $\Sigma$ the quadric:
$$\Sigma = \left\{ z_0^2 + \dots + z_{n+1}^2 = 0 \right\} \subset
{\mathbb{C}}P^{n+1}.$$
The function $\varphi_{_{\mathcal{P}}}:
{\mathbb{C}}P^{n+1} \setminus \Sigma \to \mathbb{R}$ is (up to a
constant factor):
$$\varphi_{_{\mathcal{P}}}([z_0: \ldots: z_{n+1}]) = \log
\frac{{\bigm|\sum_{j=0}^{n+1} z_j^2 \bigm|}^2}{(\sum_{j=0}^{n+1}
  |z_j|^2)^2}.$$
A straightforward computation shows that
$\varphi_{_{\mathcal{P}}}$ is Morse-Bott and
$\textnormal{Crit}(\varphi_{_{\mathcal{P}}}) = \mathbb{R}P^{n+1}$,
where $\mathbb{R}P^{n+1}$ is embedded in ${\mathbb{C}}P^{n+1}$ as
$$\mathbb{R}P^{n+1} = \left\{ [z_0: \ldots: z_{n+1}] \in
   {\mathbb{C}}P^{n+1} \bigm| z_j \in \mathbb{R} \quad \textnormal{for
     every } j \right\}.$$
Note also that $\varphi_{_{\mathcal{P}}}$
has property $(\mathcal{S})$ (see~\ref{sbsb:cond_S}).  Let $p=[1:0:
\ldots :0] \in \mathbb{R}P^{n+1}$.  A straightforward computation
shows that
$$\nabla_{\varphi}^{\textnormal{crit}} = W_p^u(X_{_{\mathcal{P}}}) =
\left\{[1:ix_1:\ldots:ix_{n+1}] \bigm| x_j \in \mathbb{R} \quad
   \textnormal{for every } j \right\} \setminus \Sigma.$$
Hence the
Lagrangian trace $\Lambda_{\mathcal{P}} \subset \Sigma$ is the
following Lagrangian sphere:
$$\Lambda_{\mathcal{P}} = \left\{[1:ix_1:\ldots:ix_{n+1}] \Bigm| x_j
   \in \mathbb{R} \quad \textnormal{for every } j, \textnormal{ and }
   \sum_{j=1}^{n+1} x_j^2 = 1 \right\}.$$

\subsubsection{Polarization of ${\mathbb{C}}P^n \times {\mathbb{C}}P^n$}
\label{sbsb:cpncpn}
Let $\mathcal{P} = (M={\mathbb{C}}P^n\times {\mathbb{C}}P^n,
\omega=\sigma \oplus \sigma, J; \Sigma)$ where
$$\Sigma = \left\{ ([z_0: \ldots :z_n],[w_0: \ldots :w_n]) \in
   {\mathbb{C}}P^n\times {\mathbb{C}}P^n \Bigm| \sum_{j=0}^{n-1} z_j
   w_j = z_n w_n \right\}.$$
A simple computation of the function
$\varphi_{_{\mathcal{P}}}$ show that it is Morse-Bott. It has one
critical submanifold which is a Lagrangian copy of ${\mathbb{C}}P^n$:
$$\Delta_{\mathcal{P}} = \left\{ ([z_0: \ldots :z_n], [\overline{z}_0:
   \ldots :\overline{z}_{n-1}: -\overline{z}_n]) \in {\mathbb{C}}P^n
   \times {\mathbb{C}}P^n \bigm| [z_0: \ldots :z_n] \in
   {\mathbb{C}}P^n \right\}.$$
Pick $p = ([0: \ldots :0:1],[0: \ldots
:0:1]) \in \Delta_{\mathcal{P}}$.  A straightforward computation shows
that
$$\nabla_{\varphi}^{\textnormal{crit}} = W_p^u(X_{_{\mathcal{P}}}) =
\left\{ ([z_0: \ldots :z_n], [\overline{z}_0: \ldots: \overline{z}_n])
   \bigm| [z_0: \ldots :z_n] \in {\mathbb{C}}P^n \textnormal{ and }
   z_n \neq 0 \right\} \setminus \Sigma.$$
Finally, the Lagrangian
trace is the following Lagrangian sphere:
$$\Lambda_{\mathcal{P}} = \left\{ ([z_0: \ldots :z_n],
   [\overline{z}_0: \ldots: \overline{z}_n]) \in \Sigma \Bigm|
   \sum_{j=0}^{n-1} |z_j|^2 = |z_n|^2 \right\}.$$

\subsubsection{Lagrangian trace of hypersurfaces in ${\mathbb{C}}P^{n+1}$}
\label{sbsb:deg-d}
Generalizing Example~\ref{sbsb:quadric} above, consider $\mathcal{P} =
(M, \omega, J; \Sigma)$ with $M={\mathbb{C}}P^{n+1}$ and $\Sigma$ the
degree $d>2$ hypersurface:
$$\Sigma = \left\{ z_0^d + \dots + z_{n+1}^d = 0 \right\} \subset
{\mathbb{C}}P^{n+1}.$$
The function $\varphi_{_{\mathcal{P}}}:
{\mathbb{C}}P^{n+1} \setminus \Sigma \to \mathbb{R}$ is (up to a
constant factor):
$$\varphi_{_{\mathcal{P}}}([z_0: \ldots: z_{n+1}]) = \log
\frac{{\bigm|\sum_{j=0}^{n+1} z_j^d \bigm|}^2}{(\sum_{j=0}^{n+1}
  |z_j|^2)^d}.$$
When $d>2$ all the critical points of
$\varphi_{_{\mathcal{P}}}$ are isolated. The critical points of index
$n+1$ are all the points $[1:\xi_1: \ldots :\xi_{n+1}]$ with
$\xi_i^d=1$ for every $i$.  Denote by $W^u_{[1: \ldots :1]}=W^u_{[1:
  \ldots :1]}(X_{_{\mathcal{P}}})$ the unstable submanifold
corresponding to the critical point $[1: \ldots :1]$, and let
$$\Lambda_{[1: \ldots :1]} = \left\{ \lim_{t \to \infty}
   Z_{_{\mathcal{P}}}^t(p) \bigm| p \in W^u_{[1: \ldots :1]} \setminus
   \{[1: \ldots :1]\} \right\}$$
be the part of the trace
corresponding to $[1: \ldots :1]$.  Here $Z_{_{\mathcal{P}}}^t$ is the
flow defined in~\ref{Sb:Lag-trace}.

We have not managed to compute $W^u_{[1: \ldots :1]}$ nor
$\Lambda_{[1: \ldots :1]}$ explicitly. However we do have the
following information on $\Lambda_{\mathcal{P}}$. Let $\xi\in
\mathbb{C}$ be a primitive root of unity of degree $d$. For every
multi-index $\underline{i} = (i_1, \ldots, i_{n+1})$ where $i_1,
\ldots, i_{n+1} \in \{ 0, \ldots, d-1\}$ denote by
$R_{\underline{i}}:{\mathbb{C}}P^{n+1} \to {\mathbb{C}}P^{n+1}$ the
map: $$R_{\underline{i}}([z_0: z_1: \ldots :z_{n+1}]) = [z_0:
\xi^{i_1} z_1: \ldots: \xi^{i_{n+1}} z_{n+1}].$$
A simple computation
shows that the vector field $X_{_{\mathcal{P}}}$ is invariant under
the action of each of the maps $R_{\underline{i}}$. Therefore the
unstable submanifold $W^u_{[1:\xi^{i_1}: \ldots: \xi^{i_{n+1}}]}$
coincides with $R_{\underline{i}}(W^u_{[1: \ldots :1]})$. We conclude
that
$$\nabla_{\varphi}^{\textnormal{crit}} = \bigcup_{\underline{i} \in I}
R_{\underline{i}}(W^u_{[1: \ldots :1]}),$$
where $I$ is the set of all
multi-indices $\underline{i} \in \{0, \ldots, d-1\}^{n+1}$.  In
particular $$\Lambda_{\mathcal{P}} = \bigcup_{\underline{i} \in I}
R_{\underline{i}}(\Lambda_{[1: \ldots :1]}),$$
hence the Lagrangian
trace $\Lambda_{\mathcal{P}}$ is a union of not more than $d^{n+1}$
possibly immersed Lagrangian spheres.

\section{Lagrangian submanifolds and polarizations} \cntrs
\label{S:Lag-pol}

In this section we consider Lagrangian submanifolds of manifolds
$\Sigma$ which appear as hyperplane sections in some polarization
$\mathcal{P} = (M, \omega, J; \Sigma)$. Given a Lagrangian $L \subset
\Sigma$ our strategy will be to go one dimension up and construct a
new Lagrangian $\Gamma_L \subset M \setminus \Sigma$. The advantage is
that, sometimes, due to the ambient geometry of $M\setminus\Sigma$ it
is easier to compute symplectic invariants of $\Gamma_L$ than those of
$L$. The basic construction is presented in detail in
subsection~\ref{Sb:Lag-circle} below. Before we continue we remark
again that all Lagrangian submanifolds are assumed to be compact and
without boundary, unless explicitly otherwise stated.

Let us briefly recall now the notions of monotone symplectic manifold
and monotone Lagrangian.  Given a symplectic manifold $(X, \omega)$ we
denote by $c_1^X \in H^2(X)$ the first Chern class of its tangent
bundle (viewed as a complex vector bundle).  A symplectic manifold
$(X, \omega)$ is called {\em spherically monotone} if the following
two conditions are satisfied:
\begin{itemize}
  \item $c_1^X$ does not vanish on $\pi_2(X)$.
  \item There exists $\lambda > 0$ such that for every $A \in
   \pi_2(X)$, $\omega(A) = \lambda c_1^X(A)$.
\end{itemize}
Denote by $N_X \in \mathbb{Z}_+$ the positive generator of the
subgroup $c_1^X(\pi_2(X)) \subset \mathbb{Z}$. We call $N_X$ the {\em
  minimal Chern number} of $(X, \omega)$.

A Lagrangian submanifold $K \subset (X, \omega)$ is called monotone if
there exists $\eta > 0$ such that the following two conditions are
satisfied:
\begin{itemize}
  \item The Maslov class of $K$, $\mu_K: \pi_2(X, K) \to \mathbb{Z}$
   is not zero.
  \item There exists $\eta > 0$ such that for every $A \in
   \pi_2(X,K)$, $\omega(A) = \eta \mu_K(A)$.
\end{itemize}
Denote by $N_K \in \mathbb{Z}_+$ the positive generator of the
subgroup $\mu_K(\pi_2(X,K)) \subset \mathbb{Z}$. We call $N_K$ the
{\em minimal Maslov number} of $K$.

\subsection{The Lagrangian circle bundle construction} \cntrsb
\label{Sb:Lag-circle}
Let $\mathcal{P} = (M, \omega, J; \Sigma)$ be a polarized K\"{a}hler
manifold. Put $\omega_{_{\Sigma}} = \omega|_{T(\Sigma)}$. Let $L
\subset (\Sigma, \omega_{_{\Sigma}})$ be a Lagrangian submanifold.
Consider the standard symplectic disc bundle $E_{\Sigma} \to \Sigma$
endowed with the symplectic structure $\frac{1}{k_{\mathcal{P}}}
\omega_{\textnormal{can}}$ (see~\ref{sbsb:sdb}). By
Theorem~\ref{T:decomposition} we can identify $(E_{\Sigma},
\frac{1}{k_{\mathcal{P}}} \omega_{\textnormal{can}})$ with $(M
\setminus \Delta_{\mathcal{P}}, \omega)$. Pick $0<\epsilon<1$, and
consider the circle bundle $P = \{ v\in E_{\Sigma} \mid
\|v\|=\epsilon\}$ over $\Sigma$. Denote by $\pi_P: P \to \Sigma$ the
projection. Finally define
$$\Gamma_L = \pi_P^{-1}(L) \subset P$$
to be the total space of the
restriction of $P$ to $L$. As we shall see in a moment
$$\Gamma_L \subset \Bigl(E_{\Sigma} \setminus \Sigma,
\frac{1}{k_{\mathcal{P}}} \omega_{\textnormal{can}} \Bigr)
\hookrightarrow (M \setminus \Sigma, \omega)$$
is a Lagrangian
submanifold.  We call $\Gamma_L$ the {\em Lagrangian circle bundle
  over $L$}.

More generally, let $\mathcal{P} = (M, \omega, J; \Sigma)$ be a
polarized K\"{a}hler manifold and $(X, \omega_{_X})$ be another
symplectic manifold. Let $L \subset (\Sigma \times X,
\omega_{_{\Sigma}} \oplus \omega_{_X})$ be a Lagrangian submanifold.
Pick $0<\epsilon<1$ and let $P = \{ v\in E_{\Sigma} \mid
\|v\|=\epsilon\}$ be as before. Consider now the circle bundle $\pi_{P
  \times X}: P \times X \to \Sigma \times X$. Define $$\Gamma_L =
\pi_{P\times X}^{-1}(L) \subset P \times X.$$
We claim that
$$\Gamma_L \subset \Bigl( (E_{\Sigma}\setminus \Sigma) \times X,
\frac{1}{k_{\mathcal{P}}} \omega_{\textnormal{can}} \oplus
(1-\epsilon^2) \omega_{_X}\Bigr) \hookrightarrow \bigl( (M \setminus
\Sigma)\times X, \omega \oplus (1-\epsilon^2)\omega_{_X} \bigr)$$
is a
Lagrangian submanifold. This follows immediately from the definition
of $\omega_{\textnormal{can}}$ (see~\ref{sbsb:sdb}) since
\begin{align*}
   & \Bigl( \frac{1}{k_{\mathcal{P}}}\omega_{\textnormal{can}}
   \oplus(1-\epsilon ^2) \omega_{_X} \Bigr){\Big|}_{T(P \times X)} =
   \bigl( \pi_P^*\omega_{_{\Sigma}} + 2r dr \wedge \alpha^{\nabla} +
   r^2 d\alpha^{\nabla} \oplus (1-\epsilon^2)\omega_{_X}
   \bigr){\big|}_{T(P \times X)} \\
   &= \bigl( (1-\epsilon^2)\pi_P^*\omega_{_{\Sigma}} \oplus
   (1-\epsilon^2)\omega_{_X} \bigl){\big|}_{T(P \times X)} =
   (1-\epsilon^2) \pi_{P\times X}^*(\omega_{_{\Sigma}} \oplus
   \omega_{_X}).
\end{align*}
Note the $(1-\epsilon^2)$ rescaling of the symplectic structure along
the $X$ factor. Note also that this construction coincides with the
previous one when $X$ is a point.  We shall call this $\Gamma_L$ too
the Lagrangian circle bundle over $L$. Note that $\Gamma_L$ depends on
the parameter $\epsilon$ (since $P$ does). In fact $\Gamma_L$'s
corresponding to different choices of $\epsilon$ are not Hamiltonianly
isotopic (even when $X=\textnormal{pt}$), however they are conformally
symplectic in $(M \setminus \Sigma) \times X$. This dependence on
$\epsilon$ will not trouble us in the sequel as $\Gamma_L$ will be
used only as an auxiliary object for studying the topology of $L$.

The following proposition compares the Maslov classes of $L\subset
\Sigma \times X$ and $\Gamma_L \subset (M \setminus \Sigma)\times X$.
We denote these classes by $\mu_L:\pi_2(\Sigma \times X, L) \to
\mathbb{Z}$, $\mu_{\Gamma_L}:\pi_2((M \setminus \Sigma) \times X,
\Gamma_L) \to \mathbb{Z}$.

\begin{prop} \label{P:Maslov}
   Let $\mathcal{P} = (M, \omega, J; \Sigma)$ be a polarized
   K\"{a}hler manifold and $(X, \omega_{_X})$ another symplectic
   manifold. Let $L \subset (\Sigma \times X, \omega_{_{\Sigma}}
   \oplus \omega_{_X})$ be Lagrangian submanifold. Let $\Gamma_L
   \subset ( (M \setminus \Sigma) \times X, \omega \oplus
   (1-\epsilon^2)\omega_{_X})$ be its corresponding circle bundle.
   Then for every $A \in \pi_2(P\times X, \Gamma_L)$,
   $$\mu_{\Gamma_L}(i_* A) = \mu_L({\pi_{P\times X}}_* A),$$
   where
   $i:P\times X \to M \times X$ is the inclusion.  Moreover, if $L
   \subset (\Sigma \times X, \omega_{_{\Sigma}} \oplus \omega_{_X})$
   is monotone and $\dim_{\mathbb{C}} \Sigma \geq 2$ then $\Gamma_L
   \subset ((M \setminus \Sigma)\times X, \omega \oplus (1-\epsilon^2)
   \omega_{_X})$ is also monotone and $N_{\Gamma_L}=N_L$.
\end{prop}

\begin{remsnonum}
   \begin{enumerate}
     \item Note that although $\Gamma_L \subset (M\setminus \Sigma)
      \times X$ is monotone (when $L$ is), $\Gamma_L$ is usually not
      monotone in $M \times X$.
     \item It can be easily seen from the proof below that statement
      on monotonicity in Proposition~\ref{P:Maslov} remains true when
      $\dim_{\mathbb{C}} \Sigma=1$ and $\mathcal{P}$ is a subcritical
      polarization. Note however that the only subcritical
      polarization $\mathcal{P} = (M, \omega, J; \Sigma)$ with
      $\dim_{\mathbb{C}}\Sigma = 1$ is $M={\mathbb{C}}P^2$,
      $\Sigma=\mathbb{C}P^1$ where $\Sigma$ is embedded as a
      projective line in ${\mathbb{C}}P^2$ (see~\cite{Bi-Ci:closed}).
   \end{enumerate}
\end{remsnonum}

\begin{proof}
   To simplify notation we present the proof for the case
   $X=\textnormal{pt}$. The proof of the general case is very similar.
   
   Let $A \in \pi_2(P, \Gamma_L)$ be represented by $\widetilde{u}:(D,
   \partial D) \to (P, \Gamma_L)$. Put $u = \pi_P \circ \widetilde{u}:
   (D, \partial D) \to (\Sigma, L)$. We have to prove that
   $\mu_{\Gamma_L}([\widetilde{u}]) = \mu_L([u])$.
   
   Put $V = M \setminus \Sigma$ and let $\xi = \ker
   (\alpha^{\nabla}|_{T(P_{\epsilon})})$ be the contact distribution
   on $P$. Denote by $N_{\Sigma} \to \Sigma$ the normal bundle of
   $\Sigma$ in $M$ (viewed as a complex line bundle). Throughout the
   proof we shall also use the following notation: given a symplectic
   vector bundle $(W, \Omega) \to X$ over a manifold $X$ and a map
   $v:Y \to X$ we shall write $(v^*W, v^* \Omega) \to Y$ for the
   pulled back symplectic vector bundle, namely for every $y \in Y$,
   $(v^*W_y, v^*\Omega_y) = (W_{v(y)}, \Omega_{v(y)})$.
   
   With these notation the symplectic vector bundle $(T(V)|_P,
   \omega)$ is isomorphic to $(\xi \oplus \pi_P^* N_{\Sigma},
   \omega|_{\xi} \oplus \pi_P^* \sigma)$, where $\sigma =
   \omega|_{N_{\Sigma}}$. Now,
   \begin{align} \label{eq:iso-vb}
      \bigl(\widetilde{u}^* T(V), \widetilde{u}^*\omega \bigr) & \cong
      \bigl(\widetilde{u}^*\xi \oplus \widetilde{u}^*\pi_P^*
      N_{\Sigma}, \widetilde{u}^*\omega|_{\xi} \oplus \widetilde{u}^*
      \pi_P^* \sigma \bigr) \cong \\
      & \cong \bigl( u^*T(\Sigma) \oplus \widetilde{u}^* \pi_P^*
      N_{\Sigma}, (1-\epsilon^2)u^* \omega_{\Sigma} \oplus
      \widetilde{u}^* \pi_P^*\sigma \bigr), \notag
   \end{align}
   where the last isomorphism follows from the fact that $\bigr(\xi,
   \omega|_{\xi}\bigl) \cong \bigl(\pi_P^* T(\Sigma),
   (1-\epsilon^2)\pi_P^*\omega_{_{\Sigma}}\bigr)$ and $\widetilde{u}^*
   \pi_P^* = u^*$.
   
   Consider now the loop of Lagrangian subspaces
   $$\widetilde{\lambda}(t) = T_{\widetilde{u}(t)}(\Gamma_L) \subset
   T_{\widetilde{u}(t)}(V), \quad t \in \partial D.$$
   In a compatible
   way to the symplectic isomorphism~\eqref{eq:iso-vb} we have the
   following isomorphism of Lagrangian subbundles over $\partial D$:
   \begin{equation} \label{eq:iso-lb}
      \widetilde{\lambda}(t) \cong T_{u(t)}(L) \oplus
      \tau_{\widetilde{u}(t)},
   \end{equation}
   where $\tau \subset T(P)$ is the subbundle $\tau = i \mathbb{R}
   \frac{\partial}{\partial r}$. The Maslov index of the loop of
   Lagrangian subspaces $\{ \tau_{\widetilde{u}(t)} \}_{t\in \partial
     D} \subset \bigl(\widetilde{u}^* \pi_P^* N_{\Sigma},
   \widetilde{u}^* \pi_P^*\sigma \bigr)$ is $0$ because the bundle
   $\tau$ is globally defined, hence this loop extends to the disc
   $D$.  Thus by~\eqref{eq:iso-lb}, $\mu_{\Gamma_L}([\widetilde{u}])$
   equals the Maslov index of the loop $\{ T_{u(t)}(L)\}_{t \in
     \partial D} \subset \bigl(u^* T(\Sigma),
   u^*\omega_{_{\Sigma}}\bigr)$ which is exactly $\mu_{L}([u])$. This
   proves the equality of the Maslov indices.
   
   Now suppose that $\dim_{\mathbb{C}} \Sigma \geq 2$. Put $n =
   \dim_\mathbb{C} M$.  As $\Delta_{\mathcal{P}} \subset M$ has
   dimension at most $n$ and $n > 2$ we have by a general position
   argument that
   $$\pi_2(M \setminus \Sigma, \Gamma_L) \cong \pi_2((M \setminus
   \Sigma) \setminus \Delta_{\mathcal{P}}, \Gamma_L) \cong \pi_2(P,
   \Gamma_L).$$
   But ${\pi_P}_*: \pi_2(P, \Gamma_L) \to \pi_2(\Sigma,
   L)$ is an isomorphism. By the equality of Maslov indices just
   proved it follows that $\Gamma_L$ is monotone if and only if $L$ is
   monotone and moreover that $N_{\Gamma_L} = N_L$.
\end{proof}

\subsection{Displaceable Lagrangian submanifolds} \cntrsb
\label{Sb:Displacable}   
Let $(\Sigma,\omega_{_{\Sigma}})$ be a closed symplectic manifold.  We
say that $(\Sigma, \omega_{_{\Sigma}})$ {\em participates in a
  polarization $\mathcal{P}$} if it can be embedded in a polarized
K\"{a}hler manifold $\mathcal{P} = (M, \omega, J; \Sigma)$ in such a
way that $\omega_{_{\Sigma}}=\omega|_{\Sigma}$.

Recall that given a polarization $\mathcal{P} = (M, \omega, J;
\Sigma)$, $(M \setminus \Sigma, J, \varphi_{_{\mathcal{P}}})$ is a
Stein manifold. Denote by $\widehat{\omega}$ the symplectic structure
associated to the completion of this Stein manifold (see
Section~\ref{Sb-Canonical} and Remark~\ref{R:completion}). Given a
Lagrangian submanifold $K$ we denote by $HF(K,K)$ the Floer homology
of $K$ with itself (see Section~\ref{S:Computations}).

\begin{thm} \label{T:Lag-subcrit}
   Let $(\Sigma, \omega_{_{\Sigma}})$ be a symplectic manifold that
   participates in a subcritical polarization $\mathcal{P} = (M,
   \omega, J; \Sigma)$. Let $(X, \omega_{_X})$ be another tame
   symplectic manifold. Let $L \subset (\Sigma \times X,
   \omega_{_{\Sigma}} \oplus \omega_{_X})$ be a Lagrangian
   submanifold. Then there exists a compactly supported $h \in
   \textnormal{Ham}((M\setminus \Sigma) \times X, \widehat{\omega}
   \oplus \omega_{_X})$ such that $h(\Gamma_L) \cap \Gamma_L =
   \emptyset$. In particular, if $L$ is monotone with $N_L \geq 2$
   then $HF(\Gamma_L, \Gamma_L) = 0$. Here the Floer homology is
   computed in $(M \setminus \Sigma) \times X$ (not in $M \times X$).
\end{thm}

\begin{remnonum}
   The Floer homology $HF(\Gamma_L, \Gamma_L)$ when computed in $(M
   \setminus \Sigma) \times X$ with either of the symplectic
   structures $\omega \oplus (1-\epsilon^2)\omega_{_X}$ or
   $\widehat{\omega} \oplus (1-\epsilon^2)\omega_{_X}$ is the same.
   This follows from convexity at infinity of $M \setminus \Sigma$.
   See~\cite{Bi-Ci:Stein} for more details on this type of argument.
\end{remnonum}

\begin{proof}[Proof of Theorem~\ref{T:Lag-subcrit}]
   Put $V = M \setminus \Sigma$, and denote by $\widehat{\omega}$ the
   symplectic structure associated to the completion of $(V, J,
   \varphi_{_{\mathcal{P}}})$.  since $\mathcal{P} = (M, \omega, J;
   \Sigma)$ is subcritical there exists a subcritical exhausting
   plurisubharmonic function $\varphi: V \to \mathbb{R}$ such that
   $(V, J, \varphi)$ is complete too (see
   Section~\ref{Sb:subcritical_polarizations} and
   Lemma~\ref{L:completion}). By Lemma~\ref{L:canonical}, $(V,
   \widehat{\omega})$ and $(V, \omega_{\varphi})$ are
   symplectomorphic.
     
   Denote by $\textnormal{pr}_{V} : V \times X \to V$ the projection.
   It now follows from Lemma~\ref{L:displacement-subcrit} that there
   exists a compactly supported $h' \in Ham(V, \widehat{\omega})$ such
   that $h'(\textnormal{pr}_V(\Gamma_L)) \cap
   \textnormal{pr}_V(\Gamma_L) = \emptyset$. Let $h = h' \times \Id_X
   : V \times X \to V \times X$.  Clearly $h(\Gamma_L) \cap \Gamma_L =
   \emptyset$. If necessary one can Hamiltonianly cut off $\Id_X$, in
   case $X$ is not compact, in order to obtain a compactly supported
   $h$.
   
   The statement on Floer homology follows immediately because by
   Proposition~\ref{P:Maslov} whenever $L$ is monotone $\Gamma_L$ is
   monotone too and $N_{\Gamma_L}=N_L$. (Note the second remark after
   the statement of Proposition~\ref{P:Maslov} for the case when
   $\dim_{\mathbb{C}} \Sigma = 1$.)
\end{proof}

\begin{thm} \label{T:Lag-crit}
   Let $(\Sigma, \omega_{_{\Sigma}})$ be a symplectic manifold that
   participates in the polarization $\mathcal{P} = (M, \omega, J;
   \Sigma)$ with property $(\mathcal{S})$ (see~\ref{Sb:Lag-trace}).
   Let $(X, \omega_{_X})$ be another tame symplectic manifold and $L
   \subset (\Sigma \times X, \omega_{_{\Sigma}} \oplus \omega_{_X})$ a
   Lagrangian submanifold.  If $L \cap (\Lambda_{\mathcal{P}} \times
   X) = \emptyset$, then there exists a compactly supported $h \in
   \textnormal{Ham}((M\setminus \Sigma) \times X, \widehat{\omega}
   \oplus (1-\epsilon^2)\omega_{_X})$ such that $h(\Gamma_L) \cap
   \Gamma_L = \emptyset$. In particular, if in addition
   $\dim_{\mathbb{C}} \Sigma \geq 2$ and $L$ is monotone with $N_L
   \geq 2$ then $HF(\Gamma_L, \Gamma_L) = 0$. Here the Floer homology
   is computed in $(M \setminus \Sigma) \times X$ (not in $M \times
   X$).
\end{thm}
\begin{remnonum}
   Theorem~\ref{T:Lag-crit} generalizes Theorem~\ref{T:Lag-subcrit}
   since if $\mathcal{P} = (M, \omega, J; \Sigma)$ is a subcritical
   polarization then $\Lambda_{\mathcal{P}} = \emptyset$.
\end{remnonum}

\begin{proof}[Proof of Theorem~\ref{T:Lag-crit}]
   The proof is similar to the one of Theorem~\ref{T:Lag-subcrit}.
   The only additional point is that, due to the assumption that $L
   \cap (\Lambda_{\mathcal{P}} \times X) = \emptyset$, we have that
   $\Gamma_L$ lies in the complement of $\nabla^{\textnormal{crit}}
   \times X$. Indeed, due to Lemma~\ref{L:completion} the Stein
   manifold $(M\setminus\Sigma, \varphi_{_{\mathcal{P}}})$ can be made
   complete in such a way that $\varphi_{_{\mathcal{P}}}$ is not
   altered in a neighbourhood of $\textnormal{pr}_M(\Gamma_L)$, where
   $\textnormal{pr}_M: M\times X \to M$ is the projection on $M$.
   Therefore, $\Gamma_L$ continues to lie in the complement of
   $\nabla^{\textnormal{crit}} \times X$ in $((M \setminus \Sigma)
   \times X, \widehat{\omega} \oplus (1-\epsilon^2) \omega_{_X})$. The
   Hamiltonian displacement along the $M \setminus \Sigma$ factor
   follows now from Lemma~\ref{L:displacement-crit}.
\end{proof}

\section{Computations in Floer homology} \cntrs
\label{S:Computations}
In this section we summarize necessary facts from Floer theory for
monotone Lagrangian submanifolds. This extension of Floer's work was
developed by Oh~\cite{Oh:HF1}. In section~\ref{Sb:spectral} we
describe a spectral sequence which converges to the Floer homology of
a Lagrangian. This spectral sequence is based on the theory developed
by Oh~\cite{Oh:spectral,Oh:relative}, however our construction is
somewhat different. We refer the reader
to~\cite{Fl:Witten-complex,Fl:Morse-theory,Oh:HF1,
  Oh:spectral,Oh:relative,Se:graded} for more details on Floer theory.

Let $(M, \omega)$ be a tame symplectic manifold (see~\cite{ALP},
chapter~10), and let $L \subset (M, \omega)$ be a monotone Lagrangian
submanifold with $N_L \geq 2$. In this situation one can define the
Floer homology of $L$ which is an invariant of the Hamiltonian isotopy
class of $L$. Let $L' = \phi(L)$, $\phi \in \textnormal{Ham}(M,
\omega)$, be a Hamiltonianly isotopic copy of $L$, and assume that $L
\pitchfork L'$.  Let $$CF(L, L') = \bigoplus_{x \in L \cap L'}
\mathbb{Z}_2 x$$
be the vector space over $\mathbb{Z}_2$ spanned by
the intersection points of $L \cap L'$. One defines a differential
$d_J: CF(L,L') \to CF(L,L')$ by choosing an almost complex structure
$J$ and counting Floer trajectories (pseudo-holomorphic strips)
connecting pairs of points of $L \cap L'$.  The homology of $d_J$,
denoted by $HF(L,L'; J)$ is called the Floer homology of the pair $(L,
L')$.

The Floer complex $CF(L,L')$ has a (relative) $\mathbb{Z} / N_L$
grading. This grading depends on a choice of a base intersection point
$x_0 \in L \cap L'$. Different choices of such a point yield a shift
in the grading. Once $x_0$ is fixed we denote by
$CF^{i(\bmod{N_L})}(L,L';x_0)$ the $i$'th $(\bmod{N_L})$ component of
$CF$.  An index computation (see~\cite{Oh:spectral}) shows that the
differential $d_J$ increases grading by $1$, $d_J: CF^{*
  (\bmod{N_L})}(L,L';x_0) \to CF^{*+1(\bmod{N_L})}(L,L';x_0)$. Thus
the Floer homology $$HF(L,L';J) = \bigoplus_{i=0}^{N_L-1}
HF^{i(\bmod{N_L})}(L,L';J,x_0)$$
has a $\mathbb{Z}/ N_L$ grading.
Again, this grading is relative as different choices of the base point
$x_0$ result in a shifted $\mathbb{Z} / N_L$ grading.  Note that there
exists a more sophisticated approach to grading, due to
Seidel~\cite{Se:graded}, which overcomes the relativity problem.

The main feature of the Floer homology is its invariance under the
choice of $L'$ (and of $J$), namely, for every $L'' = \psi(L)$, $\psi
\in \textnormal{Ham}(M,\omega)$, intersecting $L$ transversely and any
generic almost complex structures $J', J''$ there is an isomorphism
$HF(L,L';J') \cong HF(L,L''; J'')$. Moreover, this isomorphism
preserves the $\mathbb{Z}/ N_L$ grading up to a shift, namely for
given choices $x_0' \in L\cap L'$ and $x_0'' \in L \cap L''$ there
exists $s$ such that $$HF^{*(\bmod{N_L})}(L,L';J',x_0') \cong
HF^{*+s(\bmod{N_L})}(L,L'';J'',x_0'').$$

Finally, in case $L,L'$ do not intersect transversely we define
$HF(L,L') = HF(L,L'_{\epsilon})$, where $L'_{\epsilon}$ is a small
Hamiltonian perturbation of $L'$ with $L'_{\epsilon} \pitchfork L$.

\subsection{The case of $HF(L,L)$} \cntrsb
\label{Sb:HFLL}
Let $L_{\epsilon}$ be a perturbation of $L$ built in a Weinstein
neighbourhood $\mathcal{U}$ of $L$ using a $C^2$-small Hamiltonian
Morse function $f:L \to \mathbb{R}$. Assume also that $f$ has exactly
one relative minimum $x_0$. Denote by $C_f^*$ the Morse complex of
$f$.  We shall use $x_0$ as a base intersection point for the Floer
complex.  From now on we shall drop $J$ and $x_0$ from the notation of
the Floer complex and Floer homology and simply write
$CF=CF(L,L_{\epsilon})$, $d=d_J$ and
$HF^{*(\bmod{N_L})}=HF^{*(\bmod{N_L})}(L,L;J,x_0)$. It is shown
in~\cite{Oh:spectral} that
$$CF^{i(\bmod{N_L})} = \bigoplus_{j \equiv i(\bmod{N_L})} C_f^j.$$
As
$d:CF^{*(\bmod{N_L})} \to CF^{*+1(\bmod{N_L})}$ we can write
$d=\sum_{j \in \mathbb{Z}} \partial_j$ where $\partial_j$ is an
operator $\partial_j: C_f^* \to C_f^{*+1-jN_L}$. An index computation
shows that $\partial_j=0$ for every $j<0$ and, due to dimension
reasons, $\partial_j=0$ also for every $j>\nu$, where $\nu =
[\frac{\dim L+1}{N_L}]$. Thus $$d=\partial_0 + \dots +
\partial_{\nu}.$$
Roughly speaking $\partial_0$ counts the Floer
trajectories that lie in the small neighbourhood $\mathcal{U}$ of $L$,
while $\partial_1, \ldots, \partial_{\nu}$ count the ``fat''
trajectories which go out of $\mathcal{U}$. Oh proves
in~\cite{Oh:spectral} that (for suitable choices of $J$ and Riemannian
metric on $L$) the operator $\partial_0:C_f^* \to C_f^{*+1}$ can be
identified with the Morse complex differential, hence
$H^*(C_f,\partial_0) \cong H^*(L;\mathbb{Z}_2)$.  (Note however that
$\partial_j$, $j\geq 1$, are not differentials in general, namely
$\partial_j \circ \partial_j$ may not be zero.)

\subsection{A spectral sequence} \cntrsb
\label{Sb:spectral}

We shall now present a spectral sequence which enables to calculate
the Floer homology $HF(L,L)$ using the operators $\partial_1, \ldots,
\partial_{\nu}$. We continue to use the shortened notation omitting
$J$ and $x_0$ in $CF$, $d$ and $HF$.

Let $A=\mathbb{Z}_2[T,T^{-1}]$ be the algebra of Laurent polynomials
over $\mathbb{Z}_2$ in the variable $T$. We define the degree of $T$
to be $N_L$. Thus
\begin{equation*}
   A=\bigoplus_{i \in \mathbb{Z}} A^i, \quad \textnormal{where} \quad
   A^i = 
   \begin{cases}
      \mathbb{Z}_2 T^{i/N_L} & i \equiv 0(\bmod{N_L}) \\
      0 & \textnormal{othewise}
   \end{cases}
\end{equation*}
Set $\widetilde{C} = C_f \otimes A$, namely $$\widetilde{C}^l =
\bigoplus_{k \in \mathbb{Z}} C_f^{l-kN_L} \otimes A^{kN_L}, \quad
\textnormal{for every } l\in \mathbb{Z},$$
and let
$\widetilde{d}:\widetilde{C}^* \to \widetilde{C}^{*+1}$ be
$\widetilde{d} = \partial_0 \otimes 1 + \partial_1 \otimes \tau +
\dots + \partial_{\nu} \otimes \tau^{\nu}$, where $\tau^i:A^* \to
A^{*+iN_L}$ is multiplication by $T^i$.  A simple algebraic
computation shows that:
\begin{enumerate}
  \item $\widetilde{d} \circ \widetilde{d} = 0$.
  \item $H^l(\widetilde{C},\widetilde{d}) \cong HF^{l(\bmod{N_L})}$
   for every $l \in \mathbb{Z}$.
\end{enumerate}
Next we define a decreasing filtration $\cdots \subset
F^{p+1}\widetilde{C} \subset F^p \widetilde{C} \subset
F^{p-1}\widetilde{C} \subset \cdots$ on $\widetilde{C}$. For every
$p\in \mathbb{Z}$ let $A_p = \bigoplus_{k\geq p} A^{kN_L}$ be the
space of Laurent polynomials of the form $\sum_{k\geq p} \alpha_k
T^k$. Define $$F^p \widetilde{C} = \widetilde{C} \otimes A_p, \quad
\textnormal{namely } F^p \widetilde{C}^l = \bigoplus_{k\geq p}
C_f^{l-kN_L} \otimes A^{kN_L} \quad \textnormal{for every } p,l \in
\mathbb{Z}.$$
Note that since $C_f^j=0$ for every $j>\dim L$ and
$j<0$, the filtration $F^p \widetilde{C}$ is bounded. Denote by
$\{E_r^{p,q}, d_r\}$ the spectral sequence defined by this filtration.

\begin{thm} \label{T:spectral-alg}
   The spectral sequence $\{E_r^{p,q}, d_r\}$ has the following
   properties:
   \begin{enumerate}
     \item $E_0^{p,q} = C_f^{p+q-pN_L} \otimes A^{pN_L}$, $d_0 =
      \partial_0 \otimes 1$.
     \item $E_1^{p,q} = H^{p+q-pN_L}(L;\mathbb{Z}_2) \otimes
      A^{pN_L}$, $d_1=[\partial_1] \otimes \tau$, where
      $$[\partial_1]:H^{p+q-pN_L}(L; \mathbb{Z}_2) \to
      H^{p+1+q-(p+1)N_L}(L; \mathbb{Z}_2)$$
      is induced from
      $\partial_1$.
     \item For every $r\geq 1$, $E_r^{p,q}$ has the form $E_r^{p,q} =
      V_r^{p,q} \otimes A^{p N_L}$ with $d_r = \delta_r \otimes
      \tau^r$, where $V_r^{p,q}$ are vector spaces over $\mathbb{Z}_2$
      and $\delta_r$ are homomorphisms $\delta_r:V_r^{p,q} \to
      V_r^{p+r,q-r+1}$ defined for every $p,q$ and satisfy $\delta_r
      \circ \delta_r = 0$.  Moreover $$V_{r+1}^{p,q} = \frac{\ker
        (\delta_r: V_r^{p,q} \to
        V_r^{p+r,q-r+1})}{\textnormal{image\,} (\delta_r: V_r^{p-r,
          q+r-1} \to V_r^{p,q})}.$$
      (For $r=0,1$ we have
      $V_0^{p,q}=C_f^{p+q-pN_L}$, $V_1^{p,q} =
      H^{p+q-pN_L}(L;\mathbb{Z}_2)$.)
     \item $\{E_r^{p,q}, d_r\}$ collapses at the $\nu+1$ step, namely
      $d_r=0$ for every $r\geq \nu+1$ (and so
      $E_r^{p,q}=E_{\infty}^{p,q}$ for every $r\geq \nu+1$). Moreover,
      the sequence converges to $HF$, i.e. $$\bigoplus_{p+q=l}
      E_{\infty}^{p,q} \cong HF^{l(\bmod{N_L})} \quad \textnormal{for
        every } l \in \mathbb{Z}.$$
     \item For every $p \in \mathbb{Z}$, $\bigoplus_{q\in \mathbb{Z}}
      E_{\infty}^{p,q} \cong HF$.
   \end{enumerate}
\end{thm}

\begin{remsnonum}
   \begin{enumerate}
     \item Our filtration $F^p \widetilde{C}$ is different than the
      filtration used by Oh~\cite{Oh:spectral,Oh:relative}. The
      filtration used by Oh comes from the following filtration on the
      singular cohomology of $L$: $\mathcal{F}^p H^*(L;\mathbb{Z}_2) =
      \bigoplus_{0 \leq j\leq n-p} H^j(L;\mathbb{Z}_2)$. Thus the
      spectral sequence in~\cite{Oh:spectral,Oh:relative} is different
      than ours.
     \item Semi-rigorous arguments suggest that the spectral sequence
      is multiplicative in the following sense. For every $r \geq 1$
      we have a product $E_r^{p,q} \otimes E_r^{p',q'} \to
      E_r^{p+p',q+q'}$ and the differential $d_r$ satisfies Leibniz
      rule with respect to this product. Moreover the product on $E_1$
      comes from the cup product on the cohomology
      $H^*(L;\mathbb{Z}_2)$.
      
      Taking the product structure in consideration in combination
      with the other techniques of this paper should lead to new
      restrictions on the topology of Lagrangian submanifolds.
   \end{enumerate}
\end{remsnonum}

\begin{proof}[Proof of Theorem~\ref{T:spectral-alg}]
   Most of the proof is purely algebraic and follows from the basic
   construction of a spectral sequence associated to a filtered
   complex. For the readers who are not familiar with spectral
   sequences we shall now give a brief summary of this construction.
   More details can be found in any basic text on spectral sequences.
   Below we follow the conventions of~\cite{McC-Spectral}.
   
   Let $(\widetilde{C}, \widetilde{d})$ be a complex with a decreasing
   filtration $F^p \widetilde{C}$, $p\in \mathbb{Z}$. (We write the
   $\, \widetilde{}$ over $C$ and $d$ to be consistent with the
   notation of our theorem.) Assume further that the filtration is
   bounded, namely for every $l \in \mathbb{Z}$, there exist $s=s(l)$,
   $t=t(l)$ such that $F^s \widetilde{C}^l=0$ and $F^t
   \widetilde{C}^l=\widetilde{C}^l$.
   
   Define for every $p,q \in \mathbb{Z}$ and $r \geq -1$:
   \begin{itemize}
     \item $Z_r^{p,q} = \{ x \in F^p \widetilde{C}^{p+q} \mid
      \widetilde{d}x \in F^{p+r}\widetilde{C}^{p+q+1} \}$,
     \item $B_r^{p,q} = \{ x \in F^p \widetilde{C}^{p+q} \mid
      x=\widetilde{d}y, \textnormal{with } y \in F^{p-r}
      \widetilde{C}^{p+q-1} \}$,
     \item $E_r^{p,q} = Z_r^{p,q} / (Z_{r-1}^{p+1,q-1} +
      B_{r-1}^{p,q})$,
     \item $Z_{\infty}^{p,q} = \{ x \in F^p \widetilde{C}^{p+q} \mid
      \widetilde{d}x=0 \}$,
     \item $B_{\infty}^{p,q} = \{ x \in F^p \widetilde{C}^{p+q} \mid x
      = \widetilde{d}y \}$,
     \item $E_{\infty}^{p,q} = Z_{\infty}^{p,q} /
      (Z_{\infty}^{p+1,q-1} + B_{\infty}^{p,q})$.
   \end{itemize}
   A simple computation shows that $\widetilde{d}$ maps $Z_r^{p,q}$
   into $Z_r^{p+r,q-r+1}$ and moreover it descends to a homomorphism
   $d_r:E_r^{p,q} \to E_r^{p+r,q-r+1}$ such that the following diagram
   commutes:
   \[
      \begin{CD}
         Z_r^{p,q} @>{\widetilde{d}}>> Z_r^{p+r, q-r+1} \\
         @VVV  @VVV \\
         E_r^{p,q} @>{d_r}>> E_r^{p+r,q-r+1}
      \end{CD}
   \]
   Here the vertical arrows are the canonical projections. As
   $\widetilde{d} \circ \widetilde{d} = 0$ we have $d_r \circ d_r =
   0$.
   
   Denote by $F^p H(\widetilde{C}, \widetilde{d})$ the induced
   filtration on the homology of $(\widetilde{C}, \widetilde{d})$,
   namely $$F^p H^l(\widetilde{C}, \widetilde{d}) = \textnormal{Image}
   (H^l(F^p \widetilde{C}, \widetilde{d}) \to H^l (\widetilde{C},
   \widetilde{d})).$$
   Note that $F^p H(\widetilde{C}, \widetilde{d})$
   is also a bounded filtration.
   
   The fundamental features of the above construction are the
   following (see~\cite{McC-Spectral}):
   \begin{enumerate}
     \item The homology of $(E_r^{*,*},d_r)$ is isomorphic to
      $E_{r+1}^{*,*}$, namely $H(E_r^{*,*},d_r) \cong E_{r+1}^{*,*}$.
     \item $E_1^{p,q} \cong H^{p+q}(F^p
      \widetilde{C}/F^{p+1}\widetilde{C})$.
     \item $E_{\infty}^{p,q} \cong F^p
      H^{p+q}(\widetilde{C},\widetilde{d}) /
      F^{p+1}H^{p+q}(\widetilde{C}, \widetilde{d})$.
     \item For every $p,q$ there exists $r_0=r_0(p,q)$ such that
      $E_{\infty}^{p,q} = E_r^{p,q}$ for every $r \geq r_0$.
   \end{enumerate}
   If our complex $\widetilde{C}$ consists of vector spaces then
   summing up $F^p H^l(\widetilde{C},\widetilde{d})/F^{p+1}
   H^l(\widetilde{C}, \widetilde{d})$ over all $p$'s gives us an
   isomorphic copy of $H^l(\widetilde{C}, \widetilde{d})$, hence by
   $(3)$, $\bigoplus_{p+q=l} E_{\infty}^{p,q} \cong
   H^l(\widetilde{C},\widetilde{d})$.
   
   We now turn to the proof of our theorem, applying the above
   construction to our complex. A simple computation shows that
   $$Z_0^{p,q}=F^p \widetilde{C}^{p+q}, \quad
   Z_{-1}^{p+1,q-1}=F^{p+1}\widetilde{C}^{p+1}, \quad
   B_{-1}^{p,q}=0.$$
   Thus in our case $E_0^{p,q}=C_f^{p+q-pN_L}\otimes
   A^{pN_L}$, and $d_0:E_0^{p,q} \to E_0^{p,q+1}$ is just $\partial_0
   \otimes 1$. This proves statement~1.
   
   To see 2, write elements $x \in F^p \widetilde{C}^{p+q}$ as finite
   sums $x=\sum_{j\geq 0} x_{p+q-(p+j)N_L} \otimes T^{p+j}$, where
   $x_{p+q-(p+j)N_L} \in C_f^{p+q-(p+j)N_L}$. A simple computation
   shows that:
   \begin{align*}
      Z_1^{p,q} &= \Bigl\{ x=\sum_{j\geq 0} x_{p+q-(p+j)N_L} \otimes
      T^{p+j} \mid \partial_0 (x_{p+q-pN_L})=0 \Bigr\} \\
      &= Z_{\partial_0}(C_f^{p+q-pN_L}) \otimes A^{pN_L} \bigoplus
      F^{p+1} \widetilde{C}^{p+q},
   \end{align*}
   where $Z_{\partial_0}(C_f^{p+q-pN_L}) = \textnormal{Ker}
   (\partial_0:C_f^{p+q-pN_L} \to C_f^{p+q+1-pN_L})$.  Moreover we
   have:
   \begin{enumerate}
     \item $Z_0^{p+1,q-1} = F^{p+1}\widetilde{C}^{p+q}$,
     \item $B_0^{p,q} = \widetilde{d} (F^p \widetilde{C}^{p+q-1})$.
   \end{enumerate}
   It follows that $Z_0^{p+1,q-1} + B_0^{p,q} =
   \partial_0(C_f^{p+q-1-pN_L})\otimes A^{pN_L} \bigoplus
   F^{p+1}\widetilde{C}^{p+q}.$ Thus $$E_1^{p,q} =
   Z_1^{p,q}/(Z_0^{p+1,q-1} + B_0^{p,q}) =
   H^{p+q-pN_L}(L;\mathbb{Z}_2) \otimes A^{pN_L}.$$
   To compute $d_1$,
   let us describe $\widetilde{d}:Z_1^{p,q} \to Z_1^{p+1,q}$. Write an
   element $x \in Z_1^{p,q}$ as $$x=x_{p+q-pN_L} \otimes T^p +
   x_{p+q-(p+1)N_L} \otimes T^{p+1} + x',$$
   where
   $\partial_0(x_{p+q-pN_L})=0$ and $x' \in
   F^{p+2}\widetilde{C}^{p+q}$. Then $$\widetilde{d}x = \bigl(
   \partial_1(x_{p+q-pN_L}) + \partial_0(x_{p+q-(p+1)N_L}\bigr)
   \otimes T^{p+1} + \widetilde{d}x'.$$
   It follows that $d_1:
   H^{p+q-pN_L}(L;\mathbb{Z}_2) \otimes A^{pN_L} \to
   H^{p+1+q-(p+1)N_L}(L;\mathbb{Z}_2)\otimes A^{(p+1)N_L}$ has the
   form $[\partial_1]\otimes \tau$, where $[\partial_1]$ is induced
   from $\partial_1$. This completes the proof of statement~2.
   
   Statement~3 follows immediately from statement~1 by induction on
   $r$. Indeed, note that $A^{jN_L}$ is $1$-dimensional for every $j$.
   Thus the homomorphism $d_r:V_r^{p,q} \otimes A^{pN_L} \to
   V_r^{p+r,q-r+1} \otimes A^{(p+r)N_L}$ must be of the form $d_r =
   \delta_r \otimes \tau^r$ where $\delta_r$ is a homomorphism
   $V_r^{p,q} \to V_r^{p+r, q-r+1}$. As $d_r \circ d_r=0$ we also have
   $\delta_r \circ \delta_r=0$. Moreover, the homologies of $d_r$ and
   of $\delta_r$ are related by $H(V_r^{p,q}\otimes A^{pN_L},d_r) =
   H(V_r^{p,q},\delta_r) \otimes A^{pN_L}$.
   
   To prove statement~4, note that statement~1 implies that for every
   $r$, $E_r^{p,q} \neq 0$ only if $0 \leq p+q-pN_L \leq \dim L$. Let
   $r \geq \nu+1$ and $0 \leq p+q-pN_L \leq \dim L$. As $d_r:E_r^{p,q}
   \to E_r^{p+r,q-r+1}$ it is enough to show that $(p+r) +
   (q-r+1)-(p+r)N_L <0$, i.e. that $p+q-pN_L+1-rN_L < 0$. But this is
   immediate because $rN_L \geq (\nu+1)N_L > \dim L + 1$. The second
   part of statement~4 follows from the general theory of spectral
   sequences outlined at the beginning of the proof and the fact that
   $H^l(\widetilde{C},\widetilde{d}) \cong HF^{l (\bmod{N_L})}$.
   
   Finally, statement~5 follows from the following symmetries of the
   spectral sequence which are easy to check:
   \begin{enumerate}
     \item $\widetilde{d} \circ \tau = \tau \circ \widetilde{d}$.
     \item $Z_r^{p+1,q+N_L-1} = \tau(Z_r^{p,q})$,
      $B_r^{p+1,q+N_L-1}=\tau(B_r^{p,q})$ for every $p,q,r$.
     \item $E_{\infty}^{p+1,q+N_L-1} = \tau(E_{\infty}^{p,q})$ for
      every $p,q$.
   \end{enumerate}

\end{proof}

\section{Proof of the main results} \cntrs
\label{S:Proofs-main}

\begin{proof}[\textbf{Proof of Theorem~\ref{T:cpn}}]
   Since $H_1(L; \mathbb{Z})$ is $2$-torsion it is easy to see that $L
   \subset {\mathbb{C}}P^n$ must be monotone with $N_L = k(n+1)$ for
   some $k \in \mathbb{N}$.
   
   Consider the polarization $\mathcal{P} = (M={\mathbb{C}}P^{n+1},
   \sigma, J; \Sigma={\mathbb{C}}P^n)$ where $\sigma$ is the standard
   K\"{a}hler form on ${\mathbb{C}}P^{n+1}$, $J$ is the standard
   complex structure and $\Sigma \subset {\mathbb{C}}P^{n+1}$ is a
   linear hyperplane. By Example~\ref{sbsb:subcritical}, $\mathcal{P}$
   is a subcritical polarization. Put $V=M \setminus \Sigma$ and
   consider the Lagrangian circle bundle $\Gamma_L \subset V$ as
   constructed in Section~\ref{Sb:Lag-circle}. By
   Proposition~\ref{P:Maslov} $\Gamma_L$ is monotone and
   $N_{\Gamma_L}=N_L=k(n+1)$. By Theorem~\ref{T:Lag-subcrit} we have
   $HF(\Gamma_L, \Gamma_L)=0$.
   
   We now claim that $N_{\Gamma_L}=n+1$ (namely $k=1$). Indeed, if
   $k\geq 2$ then, due to dimension reasons, $d=\partial_0$, hence
   $HF(\Gamma_L, \Gamma_L) = \bigoplus_{i=0}^{n+1} H^i(\Gamma_L;
   \mathbb{Z}_2) \neq 0$. Contradiction. This proves that
   $N_{\Gamma_L}=n+1$, hence $d=\partial_0 + \partial_1$.
   
   Let $\{ E_r^{*,*}, d_r \}$ be the spectral sequence of
   Section~\ref{Sb:spectral}. By Theorem~\ref{T:spectral-alg} the
   sequence collapses at stage $r=2$, hence $E_2^{*,*} = \ldots =
   E_{\infty}^{*,*}=0$. In particular, the following sequence is
   exact: $\cdots \xrightarrow{d_1} E_1^{-1,q} \xrightarrow{d_1}
   E_1^{0,q} \xrightarrow{d_1} E_1^{1,q} \xrightarrow{d_1} \cdots$.
   Substituting $E_1^{p,q}=H^{p+q-pN_{\Gamma_L}}(\Gamma_L;
   \mathbb{Z}_2) \otimes A^{pN_{\Gamma_L}}$, $d_1=[\partial_1]\otimes
   \tau$ we obtain:
   \begin{align*}
      H^1(\Gamma_L; \mathbb{Z}_2) & \cong H^{n+1}(\Gamma_L;
      \mathbb{Z}_2) = \mathbb{Z}_2, \\
      H^n(\Gamma_L; \mathbb{Z}_2) & \cong H^0(\Gamma_L; \mathbb{Z}_2)
      =
      \mathbb{Z}_2, \\
      H^i(\Gamma_L; \mathbb{Z}_2) & = 0 \quad \textnormal{for every }
      1<i<n.
   \end{align*}
   
   In order to recover the cohomology of $L$ itself we use the Gysin
   sequence of the circle bundle $\Gamma_L \to L$. Note that the
   second Stiefel-Whitney class of the vector bundle corresponding to
   $\Gamma_L \to L$ is just $\alpha=a|_L \in H^2(L; \mathbb{Z}_2)$,
   where $a \in H^2({\mathbb{C}}P^n; \mathbb{Z}_2)$ is the generator.
   Next note that $H_1(L; \mathbb{Z})$ cannot be $0$ since if it were
   than $N_L=2(n+1)$. Thus $H_1(L; \mathbb{Z})$ is a non-trivial
   $2$-torsion group hence $H^1(L; \mathbb{Z}_2) \neq 0$.
   
   Substituting this into the $\mathbb{Z}_2$-coefficients Gysin
   sequence we obtain that $H^i(L;\mathbb{Z}_2) \xrightarrow{\cup
     \alpha} H^{i+2}(L;\mathbb{Z}_2)$ is an isomorphism for every
   $0\leq i \leq n-2$ and that $H^1(L;\mathbb{Z}_2) \cong
   H^1(\Gamma_L; \mathbb{Z}_2)$. But $H^1(\Gamma_L;
   \mathbb{Z}_2)=\mathbb{Z}_2$ hence we conclude that $H^i(L;
   \mathbb{Z}_2)=\mathbb{Z}_2$ for every $0\leq i \leq n$, which,
   additively, is precisely the $\mathbb{Z}_2$-cohomology of
   $\mathbb{R}P^n$.
   
   Finally, suppose that $n$ is even, say $n=2m$. Denote by $\beta \in
   H^1(L;\mathbb{Z}_2)$ the generator. We have seen that $\beta \cup
   \alpha^{m-1} \neq 0 \in H^{n-1}(L;\mathbb{Z}_2)$. By Poincar\'{e}
   duality, $\beta \cup \beta \cup \alpha^{m-1} \neq 0 \in
   H^n(L;\mathbb{Z}_2)$, hence $\beta \cup \beta \neq 0$. Therefore
   $\beta \cup \beta = \alpha$, and it follows that $\beta$ generates
   the cohomology ring of $L$, exactly as for $\mathbb{R}P^n$.
\end{proof}

\begin{proof}[\textbf{Proof of Theorem~\ref{T:cpnX}}]
   Consider ${\mathbb{C}}P^n \times X \subset {\mathbb{C}}P^{n+1}
   \times X$, and let $\Gamma_L \subset ({\mathbb{C}}P^{n+1} \setminus
   {\mathbb{C}}P^n) \times X$ be the Lagrangian circle bundle
   corresponding to $L \subset {\mathbb{C}}P^n \times X$.  By our
   assumptions on $L$ and on $X$, $L$ is monotone and $N_L = 2(n+1)$.
   By Proposition~\ref{P:Maslov}, $\Gamma_L$ is monotone too and
   $N_{\Gamma_L}=2(n+1)$. (See the second remark after
   Proposition~\ref{P:Maslov} for the case $n=1$.)
   
   The rest of the proof is very similar to that of
   Theorem~\ref{T:cpn}. From Theorem~\ref{T:Lag-subcrit} we obtain
   $HF(\Gamma_L, \Gamma_L)=0$. Then a similar computation via the
   spectral sequence gives:
   \begin{align*}
      H^1(\Gamma_L; \mathbb{Z}_2) & \cong H^{2n+2}(\Gamma_L;
      \mathbb{Z}_2) = \mathbb{Z}_2, \\
      H^{2n+1}(\Gamma_L; \mathbb{Z}_2) & \cong H^0(\Gamma_L;
      \mathbb{Z}_2) =
      \mathbb{Z}_2, \\
      H^i(\Gamma_L; \mathbb{Z}_2) & = 0 \quad \textnormal{for every }
      1<i<2n+1.
   \end{align*}
   The proof now continues in the same way, using the Gysin sequence,
   only that now $H^1(L;\mathbb{Z}_2)=0$. This implies that
   $H^i(L;\mathbb{Z}_2) = 0$ for every $0 < i < 2n+1$, which is
   exactly the $\mathbb{Z}_2$-cohomology of $S^{2n+1}$.
\end{proof}

\begin{proof}[\textbf{Proof of Theorem~\ref{T:cpncpn}}]
   Since $H_1(L;\mathbb{Z}_2)=0$, it is easy to see that $L \subset
   {\mathbb{C}}P^n \times {\mathbb{C}}P^n$ is monotone with
   $N_L=2(n+1)$.  Consider ${\mathbb{C}}P^n \times {\mathbb{C}}P^n
   \subset {\mathbb{C}}P^{n+1} \times {\mathbb{C}}P^n$ and let
   $\Gamma_L \subset ({\mathbb{C}}P^{n+1} \setminus {\mathbb{C}}P^n)
   \times {\mathbb{C}}P^n$ be the Lagrangian circle bundle over $L$.
   By Proposition~\ref{P:Maslov}, $\Gamma_L$ is monotone too and
   $N_{\Gamma_L}=2(n+1)$ (see the second remark after
   Proposition~\ref{P:Maslov} for the case $n=1$). By
   Theorem~\ref{T:Lag-subcrit} we have $HF(\Gamma_L, \Gamma_L)=0$.  As
   in the proof of Theorem~\ref{T:cpn}, the spectral sequence
   $\{E_r^{*,*},d_r\}$ collapses at stage $r=2$, hence $E_2^{*,*} =
   \ldots = E_{\infty}^{*,*}=0$, and we obtain the following exact
   sequences for every $q\in \mathbb{Z}$:
   $$\cdots \xrightarrow{[\partial_1]}
   H^{q-1+N_{\Gamma_L}}(\Gamma_L;\mathbb{Z}_2)
   \xrightarrow{[\partial_1]} H^q(\Gamma_L;\mathbb{Z}_2)
   \xrightarrow{[\partial_1]}
   H^{q+1-N_{\Gamma_L}}(\Gamma_L;\mathbb{Z}_2)
   \xrightarrow{[\partial_1]} \cdots$$
   As $\dim \Gamma_L = 2n+1$ and
   $N_{\Gamma_L}=2n+2$ we get $H^i(\Gamma_L;\mathbb{Z}_2)=0$ for every
   $0<i<2n+1$ and $H^0(\Gamma_L; \mathbb{Z}_2)=H^{2n+1}(\Gamma_L;
   \mathbb{Z}_2)=\mathbb{Z}_2$. Substituting this into the Gysin
   sequence of the circle bundle $\Gamma_L \to L$ we obtain that
   $H^i(L;\mathbb{Z}_2) \xrightarrow{\cup \alpha}
   H^{i+2}(L;\mathbb{Z}_2)$ is an isomorphism for every $0\leq i \leq
   2n-2$, where $\alpha \in H^2(L;\mathbb{Z}_2)$ is the second
   Stiefel-Whitney class of the vector bundle corresponding to
   $\Gamma_L \to L$. It follows that
   $H^i(L;\mathbb{Z}_2)=\mathbb{Z}_2$ for every $i=\textnormal{even}$
   and $H^i(L;\mathbb{Z}_2)=0$ for every $i=\textnormal{odd}$.
   Moreover, $\alpha \in H^2(L;\mathbb{Z}_2)$ clearly generates the
   algebra $H^*(L;\mathbb{Z}_2)$, exactly as for $H^*({\mathbb{C}}P^n;
   \mathbb{Z}_2)$.
   
   Finally, it is easy to see that $\alpha = a|_L$ where $a\in
   H^2({\mathbb{C}}P^n \times {\mathbb{C}}P^n; \mathbb{Z}_2)$ is the
   generator of $H^2$ of any of the factors of ${\mathbb{C}}P^n \times
   {\mathbb{C}}P^n$. (As $L$ is Lagrangian, it does not matter which
   factor we take.)
\end{proof}

Before we go on to the proofs of the rest of the theorems we shall
need the following proposition.
\begin{prop} \label{P:HF-spheres}
   Let $(V^{2k}, \omega)$ be a tame symplectic manifold of dimension
   $k \geq 2$ and $K^k \subset (V^{2k}, \omega)$ a monotone Lagrangian
   submanifold with $N_K \geq 2$. Suppose that $HF(K,K)=0$. Then:
   \begin{enumerate}
     \item If $K$ is a $\mathbb{Z}_2$-homology sphere, namely
      $H^*(K;\mathbb{Z}_2) \cong H^*(S^k; \mathbb{Z}_2)$, then $N_K
      \mid k+1$.
     \item If $H^1(K; \mathbb{Z}_2) \neq 0$, $H^i(K;\mathbb{Z}_2)=0$
      for every $i \neq 0,1, k-1,k$, and $k\geq 3$, $N_K\geq 3$, then
      $N_K \mid k$.
   \end{enumerate}
\end{prop}
\begin{remnonum}
   The cohomological condition in statement~2 is satisfied whenever
   $K$ is a circle bundle over a sphere of dimension $\geq 3$.
\end{remnonum}
\begin{proof}[\textbf{Proof of Proposition~\ref{P:HF-spheres}}]
   We start by proving the second statement.  Consider the spectral
   sequence $\{ E_r^{*,*}, d_r \}$ of Section~\ref{Sb:spectral}. Note
   that since $\dim K = k$, $E_r^{-r,k+r-1}=0$ for every $r\geq 1$.
   Therefore at the $r$'th step the differential $d_r$ behave as
   follows:
   \begin{equation} \label{eq:HF-spheres}
      0 \xrightarrow{d_r} E_r^{0,k} \xrightarrow{d_r} E_r^{r,k-r+1}.
   \end{equation}
   
   Since $HF(K,K)=0$, by Theorem~\ref{T:spectral-alg} (statement~5) we
   have $E_{\infty}^{*,*}=0$. On the other hand
   $E_1^{0,k}=H^k(K;\mathbb{Z}_2) \neq 0$. Denote by $r_0$ the minimal
   $r\geq 1$ for which $d_r:E_r^{0,k} \to E_r^{r,k-r+1}$ is not $0$.
   It follows from Theorem~\ref{T:spectral-alg}
   and~\eqref{eq:HF-spheres} that $E_{r_0}^{0,k} \cong H^k(K;
   \mathbb{Z}_2)$ and the homomorphism
   $$d_{r_0}=\delta_{r_0}\otimes \tau^{r_0}: H^k(K;\mathbb{Z}_2)
   \otimes A^0 \to V_{r_0}^{r_0,k-r_0+1} \otimes A^{r_0 N_K}$$
   is not
   $0$. As $H^i(K; \mathbb{Z}_2)=0$ for every $i\neq 0, 1, k-1, k$, we
   have $V_{r_0}^{r_0,k-r_0+1} \neq 0$ only if $r_0 - r_0 N_K + k
   -r_0+1 \in \{0,1,k-1,k\}$. As $N_K \geq 3$ it follows that:
   \begin{enumerate}
     \item Either $k+1 -r_0 N_K=1$, hence $N_K \mid k$;
     \item Or $k+1-r_0 N_K=0$, hence $N_K \mid k+1$.
   \end{enumerate}
   
   Since $H^1(K; \mathbb{Z}_2) \neq 0$,
   $E_1^{0,k-1}=H^{k-1}(K;\mathbb{Z}_2) \otimes A^0$ is also not $0$
   by Poincar\'{e} duality.  Applying the same arguments as above this
   time to
   $$0\xrightarrow{d_r} E_r^{0,k-1} \xrightarrow{d_r} E_r^{r,k-r}$$
   and denoting by $r_1$ the minimal $r\geq 1$ for which
   $d_r:E_r^{0,k-1} \to E_r^{r,k-r}$ is not $0$ we obtain:
   \begin{itemize}
     \item[(1')] Either $k-r_1 N_K=0$, hence $N_K \mid k$;
     \item[(2')] Or $k-r_1 N_K=1$, hence $N_K \mid k-1$.
   \end{itemize}
   Comparing cases (1),(2) with (1'),(2') above we conclude that the
   only possibility is $N_K \mid k$.
   
   The proof of statement~1 is similar (and actually simpler).
\end{proof}

\begin{proof}[\textbf{Proof of Theorem~\ref{T:cpnX-sphere}}]
   Put $m=\dim_{\mathbb{C}} X$. We shall first assume that $n+m \geq
   3$, i.e. that $\dim L \geq 3$.
   
   Embed ${\mathbb{C}}P^n \times X \subset {\mathbb{C}}P^{n+1} \times
   X$ and let $\Gamma_L \subset ({\mathbb{C}}P^{n+1} \setminus
   {\mathbb{C}}P^n) \times X$ be the Lagrangian circle bundle over
   $L$. As in the proof of Theorem~\ref{T:cpnX}, $L$ and $\Gamma_L$
   are both monotone with $N_L=N_{\Gamma_L}=2n+2$. (Note that $L$
   being a sphere of dimension $n+m \geq 3$ is simply connected.) By
   Theorem~\ref{T:Lag-subcrit}, $HF(\Gamma_L, \Gamma_L)=0$.
   
   As $\Gamma_L$ is a circle bundle over a sphere of dimension $\geq
   3$ we have $H^1(\Gamma_L; \mathbb{Z}_2) \neq 0$ and $H^i(\Gamma_L;
   \mathbb{Z}_2)=0$ for every $i\neq 0,1, n+m, n+m+1$.  Note also that
   $N_{\Gamma_L} \geq 4$. By Proposition~\ref{P:HF-spheres}, $2n+2
   \mid n+m+1$ or equivalently $m \equiv n+1(\bmod{2n+2})$.
   
   It remains to deal with the case $n+m=2$, namely $L$ being a
   Lagrangian sphere in $\mathbb{C}P^1 \times X$ where
   $\dim_{\mathbb{C}} X=1$. We claim that this is impossible under the
   assumptions of the Theorem. Indeed, $\pi_2(X)=0$ hence the homotopy
   class $[L] \in \pi_2(\mathbb{C}P^1 \times X) \cong
   \pi_2(\mathbb{C}P^1)$ comes entirely from $\pi_2(\mathbb{C}P^1)$.
   But $\omega([L])=0$ hence $[L]=0$ which is impossible since a
   Lagrangian $2$-sphere must have self-intersection $-2$.
\end{proof}

\begin{proof}[\textbf{Proof of Theorem~\ref{T:M-cover}}]
   Theorem~\ref{T:M-cover} is a special case of
   Theorem~\ref{T:M-cover-2} which will be proved in
   Section~\ref{S:Generalizations} below.
\end{proof}

\begin{proof}[\textbf{Proof of Theorem~\ref{T:cpnM-sphere}}]
   The proof is very similar to that of Theorem~\ref{T:cpnX-sphere},
   only that now $N_L = 2\gcd(n+1, N_M)$.
\end{proof}

\begin{proof}[\textbf{Proof of Theorem~\ref{T:AQ}}]
   Consider the polarization $\mathcal{P} = (M={\mathbb{C}}P^{n+1},
   \sigma, J; \Sigma=Q^n)$, where $\Sigma=Q^n$ is the quadric. By
   Example~\ref{sbsb:quadric}, the Lagrangian trace is exactly the
   Lagrangian sphere $\Lambda_Q \subset \Sigma$.
   
   Let $L \subset \Sigma$ be a Lagrangian submanifold with
   $H_1(L;\mathbb{Z})$ either $0$ or a non-trivial $2$-torsion group,
   and assume that $L \cap \Lambda_Q = \emptyset$. The minimal Chern
   number of $\Sigma$ is $n$, hence $N_L = k n$ for some $k \in
   \mathbb{N}$. Let $\Gamma_L \subset M \setminus \Sigma$ be the
   Lagrangian circle bundle over $L$. By Proposition~\ref{P:Maslov},
   $N_{\Gamma_L}=N_L=kn$. We first claim that $k=1$. Indeed, by
   Theorem~\ref{T:Lag-crit}, $HF(\Gamma_L,\Gamma_L)=0$. Hence, if
   $k\geq 2$ then since $n\geq 3$ all the differentials $d_r:
   E_r^{*,*} \to E_r^{*+r, *-r+1}$ of the spectral sequence must
   vanish for every $r \geq 1$ which is impossible since
   $E_{\infty}^{*,*}=0$ and $E_1^{*,*}$ is not $0$. This proves
   $N_{\Gamma_L}=N_L=n$. Note that this implies that
   $H_1(L;\mathbb{Z}) \neq 0$, for otherwise $N_L$ would be $2n$.
   
   By dimension reasons we have $d=\partial_0 + \partial_1$, hence
   $E_2^{*,*} = E_{\infty}^{*,*} = 0$. Following the differentials in
   the spectral sequence we obtain from Theorem~\ref{T:spectral-alg}
   the following exact sequences for every $q \in \mathbb{Z}$:
   \begin{equation} \label{eq:AQ}
      H^{q-1+n}(\Gamma_L;\mathbb{Z}_2) \xrightarrow{[\partial_1]}
      H^q(\Gamma_L;\mathbb{Z}_2) \xrightarrow{[\partial_1]}
      H^{q+1-n}(\Gamma_L;\mathbb{Z}_2).
   \end{equation}
   
   Assume first that $n>3$. From~\eqref{eq:AQ} we obtain:
   $$H^2(\Gamma_L;\mathbb{Z}_2) \cong H^{n+1}(\Gamma_L;\mathbb{Z}_2) =
   \mathbb{Z}_2, \quad H^{n-1}(\Gamma_L;\mathbb{Z}_2) \cong
   H^0(\Gamma_L;\mathbb{Z}_2)=\mathbb{Z}_2,$$
   and
   $H^i(\Gamma_L;\mathbb{Z}_2)=0$ for every $2<i<n-1$.
   
   Consider now the Gysin sequence of the circle bundle $\Gamma_L \to
   L$. The first Chern class of the normal bundle $N_{\Sigma/M}$ is
   $c_1^{N_{\Sigma/M}} = 2a|_{\Sigma}$ where $a \in
   H^2({\mathbb{C}}P^{n+1};\mathbb{Z})$ is the positive generator.
   Therefore the second Stiefel-Whitney class of the vector bundle
   associated to $\Gamma_L \to L$ is $0$. From the Gysin sequence we
   get $H^i(L;\mathbb{Z}_2)=0$ for every $2<i<n-2$, as well as the
   following exact sequence:
   $$0 \to H^2(L;\mathbb{Z}_2) \to \mathbb{Z}_2 \to
   H^1(L;\mathbb{Z}_2) \to 0.$$
   
   Since $H_1(L;\mathbb{Z}) \neq 0$ is $2$-torsion we conclude that
   $H^1(L;\mathbb{Z}_2)\neq 0$. (Actually this follows also from the
   fact that $N_L=n$ and $N_Q=n$.)  It now easily follows that
   $H^1(L;\mathbb{Z}_2)=\mathbb{Z}_2$ and that
   $H^2(L;\mathbb{Z}_2)=0$. This proves that $H^*(L;\mathbb{Z}_2)
   \cong A^*_Q$.
   
   Assume now that $n=3$. From~\eqref{eq:AQ} we obtain
   $H^2(\Gamma_L;\mathbb{Z}_2) = \mathbb{Z}_2 \oplus \mathbb{Z}_2$.
   Using the Gysin sequence and similar arguments to the preceding
   ones we obtain $H^i(L;\mathbb{Z}_2)=\mathbb{Z}_2$ for every $0\leq
   i \leq 3$, hence $H^*(L;\mathbb{Z}_2) \cong A^*_Q$.

   \begin{remsnonum}
      \begin{enumerate}
        \item Using the multiplicative structure of the spectral
         sequence (see the second remark after
         Theorem~\ref{T:spectral-alg}), it seems that the following
         should hold for $n=3$: if $a \in H^1(L;\mathbb{Z}_2)$ is the
         generator then $a \cup a = 0$.
        \item The same proof as above with small changes actually
         shows that if $L \subset Q^n$, $n\geq 3$, is a monotone
         Lagrangian submanifold with $N_L=n$ then either
         $H^*(L;\mathbb{Z}_2) \cong A_Q^*$ or $L \cap \Lambda_Q \neq
         \emptyset$.
      \end{enumerate}
   \end{remsnonum}

\end{proof}

\begin{proof}[\textbf{Proof of Theorem~\ref{T:hypersurface}}]
   Throughout the proof we set $\Sigma=\Sigma_d^n$. Consider the
   polarization $\mathcal{P} = (M={\mathbb{C}}P^{n+1}, \sigma, J;
   \Sigma)$. By Example~\ref{sbsb:deg-d} the Lagrangian trace
   $\Lambda_{\mathcal{P}} \subset \Sigma$ consists of at most
   $d^{n+1}$ immersed Lagrangian spheres.
   
   Let $L \subset \Sigma$ be a monotone Lagrangian submanifold and
   suppose that $L \cap \Lambda_{\mathcal{P}} = \emptyset$. Let
   $\Gamma_L \subset {\mathbb{C}}P^{n+1} \setminus \Sigma$ be the
   Lagrangian circle bundle over $L$. By Theorem~\ref{T:Lag-crit},
   $HF(\Gamma_L, \Gamma_L)=0$.
   
   To prove statement~1 assume first that $2d \leq n+1$ and let $L
   \subset \Sigma$ be a Lagrangian submanifold with
   $H_1(L;\mathbb{Z})=0$. A simple computation shows that
   $N_L=2(n+2-d)$ hence $N_L \geq n+3$. As $\dim \Gamma_L = n+1$ it
   follows that the Floer differential is $d=\partial_0$, hence
   $HF(\Gamma_L,\Gamma_L) \cong H^*(\Gamma_L;\mathbb{Z}_2) \neq 0$
   which is a contradiction. This proves that $L \cap
   \Lambda_{\mathcal{P}} \neq \emptyset$.  The second part of
   statement~1 (i.e. $d \geq \frac{3}{2}(n+1)$) will be treated
   towards the end of the proof.
   
   To prove statements~2 and 3, assume that $2d \leq n+1$ and that
   $H_1(L;\mathbb{Z})$ is $2$-torsion. A simple computation shows that
   $N_L=k(n+2-d)$ for some $k \in \mathbb{N}$. We first claim that
   $k=1$. Indeed, if $k\geq 2$ then $N_L \geq n+3$ and we would obtain
   $HF(\Gamma_L, \Gamma_L) \cong H^*(\Gamma_L;\mathbb{Z}_2)$
   contradicting the vanishing of Floer homology. Thus $N_L=n+2-d$. As
   $2d \leq n+1$ we have $d=\partial_0+\partial_1$ and so $E_2^{*,*} =
   E_{\infty}^{*,*}=0$. Following the differentials of the spectral
   sequence we get:
   \begin{align} \label{eq:hypersurface}
      & H^{d+1}(\Gamma_L;\mathbb{Z}_2) = \ldots =
      H^{n-d}(\Gamma_L;\mathbb{Z}_2) =0, \\
      & H^i(\Gamma_L;\mathbb{Z}_2) \cong
      H^{n-d+1+i}(\Gamma_L;\mathbb{Z}_2), \quad \textnormal{for every
      } 0 \leq i \leq d. \notag
   \end{align}
 
   Consider now the Gysin sequence of $\Gamma_L \to L$. Note that the
   second Stiefel-Whitney class $w$ of the vector bundle corresponding
   to $\Gamma_L \to L$ is $d a|_L \in H^2(L;\mathbb{Z}_2)$ where $a
   \in H^2({\mathbb{C}}P^{n+1};\mathbb{Z}_2)$ is the generator.  When
   $d=$ even we have $w=0$, hence the Gysin sequence gives the exact
   sequence:
   $$0 \to H^i(L;\mathbb{Z}_2) \to H^i(\Gamma_L;\mathbb{Z}_2) \to
   H^{i-1}(L;\mathbb{Z}_2) \to 0 \quad \textnormal{for every } i.$$
   Combining with~\eqref{eq:hypersurface} we obtain
   $H^j(L;\mathbb{Z}_2)=0$ for every $d \leq j \leq n-d$. Moreover, we
   have $\beta_i(L)+\beta_{i-1}(L)=\beta_i(\Gamma_L)$.
   By~\eqref{eq:hypersurface} and Poincar\'{e} duality we have
   $\beta_{d-i}(\Gamma_L)=\beta_i(\Gamma_L)$ for every $0\leq i \leq
   d$, hence $$\beta_i(L) + \beta_{i-1}(L) =
   \beta_{d-i}(L)+\beta_{d-i-1}(L) \quad \textnormal{for every } 0\leq
   i \leq d.$$
   Putting $i=0$ we obtain $\beta_0(L)=\beta_{d-1}(L)$
   (because $\beta_d(L)=0$). Next,
   $\beta_1(L)+\beta_0(L)=\beta_{d-1}(L)+\beta_{d-2}(L)$, hence
   $\beta_1(L)=\beta_{d-2}(L)$. Continuing by induction we obtain
   $\beta_i(L) = \beta_{d-1-i}(L)$ for every $0\leq i \leq d-1$.
   
   The statement for $d=$ odd follows at once
   from~\eqref{eq:hypersurface} and the fact that the Stiefel-Whitney
   class $w = d a|_L$ equals $a|_L$ since $d$ is odd.
   
   Statement~4 follows easily from statement~2 of
   Proposition~\ref{P:HF-spheres} applied to $\Gamma_L$.
   
   To prove the second part of statement~1 and statement~5 we shall
   use an extension of Floer homology for so called {\em strongly
     negative} Lagrangian submanifolds, due to
   Lazzarini~\cite{Laz:discs}. Let $K \subset (V, \omega)$ be a
   Lagrangian in a tame symplectic manifold. $K$ is called strongly
   negative if there exists $\lambda < 0$ such that $\mu_K = \lambda
   [\omega]$ on $\pi_2(V,K)$ and in addition the following conditions
   are satisfied:
   \begin{enumerate}
     \item $c_1^V(A) \leq 2 - \dim_{\mathbb{C}}V$ for every $A \in
      \pi_2(V)$ with $\omega(A)>0$.
     \item $\mu_K(A) \leq 2 - \dim_{\mathbb{C}}V$ for every $A \in
      \pi_2(V,K)$ with $\omega(A)>0$.
   \end{enumerate}
   Under these assumptions the Floer homology $HF(K,K')$ is well
   defined for every Lagrangian $K'$ which is Hamiltonianly isotopic
   to $K$ and moreover $$HF(K,K) \cong HF(K,K') \cong
   H^*(K;\mathbb{Z}_2).$$
   The reason for this, roughly speaking, is
   that under the above negativity assumptions there are no
   $J$-holomorphic spheres or discs with boundary on $K$ for generic
   almost complex structure $J$. Hence, the Floer differential for
   $CF(K,K)$ is the Morse-homology differential $d=\partial_0$. The
   non-existence of spheres and discs is due to negative dimension of
   the moduli spaces of $J$-holomorphic spheres and discs. For the
   dimension formulae to hold one has to work with regular almost
   complex structures $J$. Regularity may be achieved by generic
   perturbations of $J$ as long as the $J$-holomorphic discs/spheres
   in question are somewhere injective (see~\cite{McD-Sa:Jhol-2}).
   Thus an essential ingredient in applying the dimension argument is
   a procedure which enables to extract a somewhere injective disc
   from a given pseudo-holomorphic disc.  Such procedures have been
   developed by Lazzarini~\cite{Laz:discs} and by Kwon and
   Oh~\cite{Kw-Oh:discs}.
   
   Coming back to the proof of the second part of statement~1, assume
   that $d \geq \frac{3}{2}(n+1)$, and $n\geq 3$. Note that since
   $H_1(L;\mathbb{Z})=0$ we have $\mu_L =
   2(n+2-d)[\omega_{_{\Sigma}}]$. As $n+2-d < 0$, for every $A \in
   \pi_2(\Sigma,L)$ with $\omega_{_{\Sigma}}(A)>0$ we have $\mu_L(A)
   \leq 2(n+2-d)$. By Proposition~\ref{P:Maslov} for every $A \in
   \pi_2({\mathbb{C}}P^{n+1} \setminus \Sigma, \Gamma_L)$ with
   $\omega(A)>0$ we have $\mu_{\Gamma_L}(A) \leq 2(n+2-d)$. But since
   $d \geq \frac{3}{2}(n+1)$ we have $$2(n+2-d) \leq 2-(n+1) = 2-
   \dim_{\mathbb{C}}({\mathbb{C}}P^{n+1}\setminus \Sigma).$$
   Thus
   $\Gamma_L \subset {\mathbb{C}}P^{n+1} \setminus \Sigma$ is strongly
   negative. It follows that $HF(\Gamma_L,\Gamma_L)\cong
   H^*(\Gamma_L;\mathbb{Z}_2)$ which contradicts the vanishing of
   Floer homology. This proves that $L \cap \Lambda_d \neq \emptyset$.
   
   We turn to the proof of statement~5. Since $H_1(L;\mathbb{Z})$ is
   $t$-torsion and $n+2-d < 0$ we have $\mu_L(A)\leq
   \frac{2(n+2-d)}{t}$ for every $A \in \pi_2(\Sigma,L)$ with
   $\omega_{_{\Sigma}}(A)>0$.  By our assumption that $d \geq
   \frac{t(n-1)}{2} + n+2$ and Proposition~\ref{P:Maslov}, for every
   $A \in \pi_2({\mathbb{C}}P^{n+1} \setminus \Sigma, \Gamma_L)$ with
   $\omega(A)>0$ we have $\mu_{\Gamma_L}(A) \leq \frac{2(n+2-d)}{t}
   \leq 1-n = 2-\dim_{\mathbb{C}}({\mathbb{C}}P^{n+1} \setminus
   \Sigma)$, hence $\Gamma_L$ is strongly negative. The rest of the
   proof is very similar to the preceding arguments.
\end{proof}

\begin{proof}[\textbf{Proof of Theorem~\ref{T:sigmacpncpn}}]
   Consider the polarization $\mathcal{P} = (M={\mathbb{C}}P^n\times
   {\mathbb{C}}P^n, \omega=\sigma \oplus \sigma, J; \Sigma)$.  By
   Example~\ref{sbsb:cpncpn} the Lagrangian trace
   $\Lambda_{\mathcal{P}}$ is exactly the Lagrangian sphere
   $\Lambda_{\Sigma}$.
   
   Let $L \subset \Sigma$ be a Lagrangian submanifold with
   $H_1(L;\mathbb{Z})=0$ and suppose that $L \cap \Lambda_{\Sigma} =
   \emptyset$. A simple computation shows that $L$ is monotone with
   $N_L=2n$. Consider the Lagrangian circle bundle $\Gamma_L \to L$.
   By Proposition~\ref{P:Maslov}, $\Gamma_L \subset
   ({\mathbb{C}}P^n\times {\mathbb{C}}P^n) \setminus \Sigma$ is
   monotone too and $N_{\Gamma_L}=N_L=2n$. By
   Theorem~\ref{T:Lag-crit}, $HF(\Gamma_L,\Gamma_L)=0$.
   
   In a similar manner to the proofs of Theorems~\ref{T:cpn}
   and~\ref{T:cpnX} we recover the cohomology of $\Gamma_L$ via the
   spectral sequence, obtaining:
   \begin{align*}
      & H^0(\Gamma_L;\mathbb{Z}_2) \cong H^1(\Gamma_L;\mathbb{Z}_2)
      \cong H^{2n-1}(\Gamma_L;\mathbb{Z}_2) \cong
      H^{2n}(\Gamma_L;\mathbb{Z}_2)
      \cong \mathbb{Z}_2, \\
      & H^i(\Gamma_L;\mathbb{Z}_2)=0 \quad \textnormal{for every }
      1<i<2n-1.
   \end{align*}
   Next note that the restriction of the complex line bundle
   $N_{\Sigma/{\mathbb{C}}P^n \times {\mathbb{C}}P^n}$ to $L$ is
   trivial, since its first Chern class is just $[\omega]|_L$ and
   $H^2(L;\mathbb{Z})$ is torsion-free because $H_1(L;\mathbb{Z})=0$.
   Thus the circle bundle $\Gamma_L \to L$ is trivial, hence
   $H^*(\Gamma_L) = H^*(L) \otimes H^*(S^1)$, from which it follows
   that $H^*(L;\mathbb{Z}_2) \cong H^*(S^{2n-2};\mathbb{Z}_2)$.
\end{proof}

\section{Generalizations} \cntrs
\label{S:Generalizations}
Below we present miscellaneous generalizations of some of the theorems
from Section~\ref{S:Intro}. Since the proofs are rather analogous to
those of Section~\ref{S:Proofs-main} we shall only outline the proofs
omitting repeated details.

The following theorem generalizes Theorem~\ref{T:M-cover}.
\begin{thm} \label{T:M-cover-2}
   Let $X$ be a symplectic manifold that has a covering which is
   symplectomorphic to a domain in a subcritical Stein manifold. Let
   $M$ be a closed spherically monotone symplectic manifold. Assume
   $\dim M, \dim X > 0$. If $M \times X$ has a Lagrangian sphere then
   $2 N_M \mid \dim_{\mathbb{C}}M + \dim_{\mathbb{C}}X + 1$.
\end{thm}
Examples of subcritical Stein manifolds, other than $\mathbb{C}^n$,
can be found in~\cite{Bi:barriers, Bi-Ci:closed}. (Note however, that
after completion all subcritical Stein manifolds are
split~\cite{Ci:split}.)

\begin{proof}[\textbf{Proof of Theorem~\ref{T:M-cover-2}}]
   Let $Y \to X$ be a covering of $X$ by a symplectic manifold $Y$
   which is (symplectomorphic to) a domain in a subcritical Stein
   manifold $V$. Let $L \subset M \times X$ be a Lagrangian sphere. By
   assumption $\dim L \geq 2$, hence $L$ is simply connected. Consider
   the lift $\widetilde{L} \subset M \times Y \subset M \times V$ of
   $L$. Clearly $\widetilde{L}$ is also an embedded Lagrangian sphere,
   and $N_{\widetilde{L}}=2N_M$. As $V$ is subcritical we have
   $HF(\widetilde{L}, \widetilde{L})=0$. By
   Proposition~\ref{P:HF-spheres}, $2N_M \mid
   \dim_{\mathbb{C}}M+\dim_\mathbb{C}X+1$.
\end{proof}

A symplectic manifold $(\Sigma, \omega_{_{\Sigma}})$ is called {\em
  monotone} if there exists $\lambda>0$ such that
$[\omega_{_{\Sigma}}] = \lambda c_1^{\Sigma} \in
H^2(\Sigma;\mathbb{R})$, where $c_1^{\Sigma}$ is the first Chern class
of the tangent bundle of $\Sigma$. We denote by $N_{\Sigma}^H \in
\mathbb{N}$ the positive generator of the subgroup $\{ c_1^{\Sigma}(A)
\mid A \in H_2(\Sigma;\mathbb{Z}) \}$.
\begin{thm} \label{T:hypersurface-2}
   Let $\Sigma$ be a closed monotone symplectic manifold (resp.
   spherically monotone) that participates in some polarization
   $\mathcal{P}$ (see Section~\ref{Sb:Displacable}). Denote
   $n=\dim_{\mathbb{C}}\Sigma$.
   \begin{enumerate}
      
     \item Suppose that $\mathcal{P}$ is subcritical. Then:
      \begin{enumerate}
        \item If $2N_{\Sigma}^H > n+1$ (resp.  $2N_{\Sigma} > n+1$)
         then there exist no Lagrangian submanifolds $L \subset
         \Sigma$ with $H_1(L;\mathbb{Z})=0$ (resp. $\pi_1(L)=0$).
        \item If $2N_{\Sigma}^H=n+1$ (resp.  $2N_{\Sigma}=n+1$) then
         every Lagrangian submanifold $L \subset \Sigma$ with
         $H_1(L;\mathbb{Z})=0$ (resp. $\pi_1(L)=0$) must satisfy
         $H^*(L;\mathbb{Z}_2) \cong H^*(S^n;\mathbb{Z}_2)$.
      \end{enumerate}
      
     \item Suppose $\mathcal{P}$ has property $(\mathcal{S})$ (see
      Section~\ref{Sb:Lag-trace}). Then:
      \begin{enumerate}
        \item If $2N_{\Sigma}^H > n+1$ (resp.  $2N_{\Sigma} > n+1$)
         then for every Lagrangian submanifold $L \subset \Sigma$ with
         $H_1(L;\mathbb{Z})=0$ (resp. $\pi_1(L)=0$) must satisfy $L
         \cap \Lambda_{\mathcal{P}} \neq \emptyset$.
        \item If $2N_{\Sigma}^H=n+1$ (resp.  $2N_{\Sigma}=n+1$) then
         every Lagrangian submanifold $L \subset \Sigma$ with
         $H_1(L;\mathbb{Z})=0$ (resp. $\pi_1(L)=0$) and with $L \cap
         \Lambda_{\mathcal{P}} = \emptyset$ must satisfy
         $H^*(L;\mathbb{Z}_2) \cong H^*(S^n;\mathbb{Z}_2)$.
      \end{enumerate}

   \end{enumerate}
\end{thm}
We remark that a similar statement to 1(a) above was proved
in~\cite{Bi-Ci:closed} using a different approach.

\begin{proof}[\textbf{Proof of Theorem~\ref{T:hypersurface-2}}]
   Note that statement~(1) is a special case of statement~(2), since
   for a subcritical polarization $\mathcal{P}$ we have
   $\Lambda_{\mathcal{P}}=\emptyset$. We therefore prove
   statement~(2).
   
   Let $L \subset \Sigma$ be a Lagrangian submanifold with
   $H_1(L;\mathbb{Z})=0$ and $L \cap \Lambda_{\mathcal{P}} =
   \emptyset$. Note that $L$ is monotone and that $N_L \geq
   2N_{\Sigma}^H$ (here we continue to denote by $N_L$ the minimal
   Maslov number, i.e. the positive generator of the subgroup
   $\mu_L(\pi_2(\Sigma,L)) \subset \mathbb{Z}$.)
   
   Consider now the Lagrangian circle bundle $\Gamma_L \to L$ in $M
   \setminus \Sigma$ (where, $\mathcal{P} = (M, \omega, J; \Sigma)$ is
   the polarization in which $\Sigma$ participates). By
   Theorem~\ref{T:Lag-crit}, $HF(\Gamma_L, \Gamma_L)=0$.
   
   As in the end of the proof of Theorem~\ref{T:sigmacpncpn}, since
   $H_1(L;\mathbb{Z})=0$ and $L \subset \Sigma$ is Lagrangian, the
   bundle $\Gamma_L \to L$ must be trivial. Thus $H^1(\Gamma_L;
   \mathbb{Z}_2) = \mathbb{Z}_2$. In view of this, the only way the
   spectral sequence can converge to $0$ is if $N_L \leq n+1$.  Thus
   if $2N_{\Sigma}^H > n+1$ we arrive at a contradiction. This
   proves~2(a).
   
   Suppose now that $2N_{\Sigma}^H=n+1$. Computing using the spectral
   sequence we obtain:
   \begin{align*}
      & H^0(\Gamma_L;\mathbb{Z}_2) \cong H^1(\Gamma_L;\mathbb{Z}_2)
      \cong H^n(\Gamma_L;\mathbb{Z}_2) \cong
      H^{n+1}(\Gamma_L;\mathbb{Z}_2) = \mathbb{Z}_2, \\
      & H^i(\Gamma_L;\mathbb{Z}_2)=0, \quad \textnormal{for every }
      1<i<n.
   \end{align*}
   It now easily follows that $H^*(L;\mathbb{Z}_2) \cong
   H^*(S^n;\mathbb{Z}_2)$.
   
   We omit the proof of the statements assuming spherical monotonicity
   since it is completely analogous to the proof above.
\end{proof}

The following theorem generalizes statement~4 of
Theorem~\ref{T:hypersurface}.
\begin{thm} \label{T:hypersurface-sphere}
   Let $\Sigma$ be a closed spherically monotone symplectic manifold
   that participates in a polarization $\mathcal{P}$ which has
   property $(\mathcal{S})$ (see Sections~\ref{Sb:Displacable}
   and~\ref{Sb:Lag-trace}). Assume that $\dim_{\mathbb{C}}\Sigma \geq
   3$ and $N_{\Sigma} \geq 2$. If $2N_{\Sigma} \nmid
   \dim_{\mathbb{C}}\Sigma+1$ then every Lagrangian sphere $L \subset
   \Sigma$ must satisfy $L \cap \Lambda_{\mathcal{P}} \neq \emptyset$.
\end{thm}
We omit the proof as it is rather similar to that of
Theorem~\ref{T:hypersurface}.

The next theorem provides ``Euler characteristic'' type restrictions
on monotone Lagrangians. Let $(V, \omega)$ be a symplectic manifold
and $L \subset (V,\omega)$ a monotone Lagrangian submanifold with
minimal Maslov number $N_L \geq 2$. For every $j \in \mathbb{Z}$,
denote by $\gamma_j$ the sum of all $\mathbb{Z}_2$-Betti numbers of
indices that are congruent to $j$ modulo $N_L$, namely
$$\gamma_j = \sum_{k \in \mathbb{Z}} \dim_{\mathbb{Z}_2}
H^{kN_L+j}(L;\mathbb{Z}_2).$$
Next, for every two integers $s\leq t$
denote:
$$\chi_{s,t}(L) = \gamma_s - \gamma_{s+1} + \cdots + (-1)^{t-s}
\gamma_t.$$
Given a Morse function $f:L \to \mathbb{R}$ and $s\in
\mathbb{Z}$, put
$$\kappa_s(f) = \# \{ p \in \textnormal{Crit}(f) \mid
\textnormal{ind}_p(f) \equiv s (\bmod {N_L}) \}.$$
Next, for for every
$s \leq t \in \mathbb{Z}$ put
\begin{align*}
   \lambda_s(L) &= \min_{f \in \textnormal{Morse}} \Bigl\{ \min \{
   \kappa_{s-1}(f), \kappa_s(f) \} \Bigr\}, \\
   \lambda_{s,t}(L) &= \min_{f \in \textnormal{Morse}} \Bigl\{ \min \{
   \kappa_{s-1}(f), \kappa_s(f)\} + \min \{ \kappa_t(f),
   \kappa_{t+1}(f) \} \Bigr\}.
\end{align*}

\begin{thm} \label{T:Euler-char}
   Let $(V,\omega)$ be a tame symplectic manifold and $L \subset (V,
   \omega)$ be a monotone Lagrangian submanifold with $N_L \geq 2$.
   Put $\nu = [\frac{\dim L + 1}{N_L}]$. Suppose $HF(L,L)=0$. Then for
   every $s \leq t$ we have:
   \begin{enumerate}
     \item If $t-s=$ even:
      \begin{enumerate}
        \item $0 \leq \chi_{s,t}(L) \leq \nu \min \{ \gamma_{s-1},
         \gamma_s \} + \nu \min \{ \gamma_t, \gamma_{t+1} \}$.
        \item $\chi_{s,t}(L) \leq \lambda_{s,t}(L)$.
      \end{enumerate}
     \item If $t-s=$ odd:
      \begin{enumerate}
        \item $-\nu \min \{\gamma_t, \gamma_{t+1} \} \leq
         \chi_{s,t}(L) \leq \nu \min \{\gamma_{s-1}, \gamma_s\}$.
        \item $-\lambda_t(L) \leq \chi_{s,t}(L) \leq \lambda_s(L)$.
      \end{enumerate}
   \end{enumerate}
   In particular, if $(V,J,\varphi)$ is a Stein manifold which is
   either subcritical, or has property $(\mathcal{S})$ of
   Section~\ref{sbsb:cond_S} and $L \cap
   \nabla_{\varphi}^{\textnormal{crit}} = \emptyset$, then the above
   inequalities hold.
\end{thm}

To prove Theorem~\ref{T:Euler-char} we shall need the following simple
Lemma from linear algebra.
\begin{lem} \label{L:Euler-char}
   Let $(D=\bigoplus_{i \in \mathbb{Z}} D^i, \partial)$ be a complex
   of vector spaces and $H(D,\partial) = \bigoplus_{i \in \mathbb{Z}}
   H^i(D, \partial)$ its cohomology. For every two integers $s\leq t$
   put
   $$
   \chi_{s,t}(D) = \sum_{i=s}^t (-1)^{i-s} \dim D^i, \quad
   \chi_{s,t}(H(D,\partial)) = \sum_{i=s}^t (-1)^{i-s} \dim
   H^i(D,\partial).$$
   Then $\chi_{s,t}(D) = \chi_{s,t}(H(D,\partial))
   + \dim \partial(D^{s-1}) + (-1)^{t-s}\dim \partial(D^t)$.
   
   In particular, for $t-s=$ even,
   $$\chi_{s,t}(H(D,\partial)) \leq \chi_{s,t}(D) \leq
   \chi_{s,t}(H(D,\partial)) + \min \{\dim D^{s-1}, \dim D^s\} + \min
   \{ \dim D^t, \dim D^{t+1}\},$$
   while for $t-s=$ odd,
   $$\chi_{s,t}(H(D,\partial)) - \min \{\dim D^t, D^{t+1}\} \leq
   \chi_{s,t}(D) \leq \chi_{s,t}(H(D,\partial)) + \min \{\dim D^{s-1},
   D^s\}.$$
\end{lem}
The proof of the lemma is completely straightforward, we therefore
omit it.

\begin{proof}[Proof of Theorem~\ref{T:Euler-char}]
   We first prove the second and fourth inequalities.  Given a Morse
   function $f:L \to \mathbb{R}$, denote by $C_f^*$ the Morse complex
   associated to $f$. For every $i\in \mathbb{Z}$, put $C^i = C_f^{i
     (\bmod{N_L})} = \bigoplus_{k \in \mathbb{Z}} C_f^{i+kN_L}$.
   Recall from Section~\ref{Sb:HFLL} that $C^*$ can be endowed with
   two differentials: the Morse differential $\partial_0$ and the
   Floer differential $d=\partial_0 + \cdots + \partial_{\nu}$. Thus
   we have $H^*(C,\partial_0) = H^{*(\bmod{N_L})}(L;\mathbb{Z}_2)$ and
   $H^*(C,d) = HF^{*(\bmod{N_L})}(L,L)$.  By Lemma~\ref{L:Euler-char}
   (applied for $\partial_0$ and for $d$), for every $s \leq t$ we
   have:
   \begin{align*}
      \chi_{s,t}(C) &= \chi_{s,t}(H(C,\partial_0)) + \dim
      \partial_0(C^{s-1}) + (-1)^{t-s}\dim \partial_0(C^t), \\
      \chi_{s,t}(C) &= \chi_{s,t}(H(C,d)) + \dim d(C^{s-1}) +
      (-1)^{t-s}d(C^t).
   \end{align*}
   By assumption $H(C,d)=HF(L,L)=0$, hence
   \begin{align} \label{eq:chi-1}
      \chi_{s,t}(L) &= \chi_{s,t}(H(C,\partial_0)) \\
      &= \dim d(C^{s-1}) - \dim \partial_0(C^{s-1}) + (-1)^{t-s}(\dim
      d(C^t) - \dim \partial_0(C^t)). \notag
   \end{align}
   Now if $t-s=$ even we get $$\chi_{s,t}(L) \leq \dim d(C^{s-1}) +
   \dim d(C^t) \leq \min \{ \kappa_{s-1}(f), \kappa_s(f) \} + \min
   \{\kappa_t(f),\kappa_{t+1}(f) \}.$$
   Taking the minimum over all
   Morse functions $f:L \to \mathbb{R}$ we obtain $\chi_{s,t}(L) \leq
   \lambda_{s,t}(L)$. Assume now that $t-s=$ odd. As $d=\partial_0 +
   \cdots + \partial_{\nu}$, a simple dimension computation (using the
   grading of each $\partial_k$) shows that $\dim d(C^i) \geq \dim
   \partial_0(C^i)$ for every $i \in \mathbb{Z}$. Using this
   with~\eqref{eq:chi-1} we get $$-\min \{\kappa_t(f),
   \kappa_{t+1}(f)\} \leq -\dim d(C^t) \leq \chi_{s,t}(L) \leq \dim
   d(C^{s-1}) \leq \min\{ \kappa_{s-1}(f), \kappa_s(f)\}.$$
   Since this
   is true for all Morse functions $f:L \to \mathbb{R}$ we obtain
   $-\lambda_t(L) \leq \chi_{s,t}(L) \leq \lambda_s(L)$.
   
   We now turn to the proof of the first and third inequalities. Let
   $\{E_r^{*,*}, d_r\}$ be the spectral sequence of
   Section~\ref{Sb:spectral}. For every $r\geq 0$, $l\in \mathbb{Z}$,
   put $\bar{E}_r^l = \bigoplus_{p+q=l} E_r^{p,q}$. Note that
   $\bar{E}^*_{r+1} = H^*(\bar{E}_r,d_r)$, and $\bar{E}_1^l \cong
   H^{l(\bmod{N_L})}(L;\mathbb{Z}_2)$ for every $l \in \mathbb{Z}$.
   By Lemma~\ref{L:Euler-char} applied $\nu$ times we obtain:
   \begin{align} \label{eq:chi-2}
      \chi_{s,t}(L)=\chi_{s,t}(\bar{E}_1) &=
      \chi_{s,t}(\bar{E}_2)+\dim
      d_1(\bar{E}_1^{s-1}) + (-1)^{t-s}\dim d_1(\bar{E}_1^t) \\
      &= \ldots = \chi_{s,t}(\bar{E}_{\nu+1}) + \sum_{r=1}^{\nu}
      \bigl(\dim d_r(\bar{E}_r^{s-1}) + (-1)^{t-s}\dim
      d_r(\bar{E}_r^t) \bigr). \notag
   \end{align}
   Now $\bar{E}_{\nu+1}=0$ because $HF(L,L)=0$. Note that by
   Theorem~\ref{T:spectral-alg} we have that $\dim d_r(\bar{E}_r^l)
   \leq \min \{\gamma_l, \gamma_{l+1}\}$ for every $l\in \mathbb{Z}$.
   The desired inequalities now easily follow from~\eqref{eq:chi-2}.
\end{proof}


\end{document}